\documentclass[12pt]{elsarticle}
\usepackage{graphicx}
\usepackage{amssymb}    %American Mathematical Society symbol package
\usepackage{epstopdf}
\usepackage{comment}
\usepackage{xcolor}
\usepackage{hyperref}
\usepackage[final]{pdfpages}
\usepackage{frcursive}
\usepackage{bbm}
\usepackage{mathtools}
\usepackage{tabularx}
\usepackage{algorithm}
\usepackage{algorithmic}
\usepackage[normalem]{ulem}

\usepackage{upgreek}

\usepackage{yfonts}
\usepackage{lscape}
\usepackage{epsfig}
\usepackage{fullpage}
\usepackage[titletoc]{appendix}
\usepackage{fancyhdr}
\usepackage{amsmath,amssymb,xspace}
\usepackage{tabularx}
\usepackage{mathrsfs}
\def\cal{\mathcal}

\newcommand{\bs}[1]{\boldsymbol{#1}}

\newcommand{\R}{\mathbb{R}}

\newcommand{\bl}{\begin{law}}
\newcommand{\el}{\end{law}}
\newcommand{\bthm}{\begin{thm}}
\newcommand{\ethm}{\end{thm}}

\newfont{\bbb}{msbm10 scaled\magstep1}

\let \leq \leqslant
\let \geq \geqslant
\let \cal \mathcal

{ \par \medskip \par
  \noindent \textit{\textbf{Demonstration\/}} : }{\null \hfill $\Box$ \par }

\makeatletter
\newcommand{\doublewidetilde}[1]{{%
  \mathpalette\double@widetilde{#1}%
}}
\newcommand{\double@widetilde}[2]{%
  \sbox\z@{$\m@th#1\widetilde{#2}$}%
  \ht\z@=.9\ht\z@
  \widetilde{\box\z@}%
}

\begin{document}

\begin{frontmatter}

\title{Non-diffusive neural network method for hyperbolic conservation laws}

\author[carl,crm]{Emmanuel LORIN}
\ead{elorin@math.carleton.ca}

\address[carl]{School of Mathematics and Statistics, Carleton University, Ottawa, Canada, K1S 5B6}
\address[crm]{Centre de Recherches Math\'{e}matiques, Universit\'{e} de Montr\'{e}al, Montr\'{e}al, Canada, H3T~1J4}

\author[ott]{Arian NOVRUZI\corref{cor1}}
\ead{novruzi@uottawa.ca}
\address[ott]{Department of Mathematics and Statistics, University of Ottawa, 
	%STEM Complex, 150 Louis-Pasteur Pvt, 
	Ottawa, ON K1N 6N5, Canada}
\cortext[cor1]{Corresponding author}

\begin{abstract}

In this paper we develop a non-diffusive neural network (NDNN) algorithm for accurately solving weak solutions to hyperbolic conservation laws. The principle is to construct these weak solutions by computing smooth local solutions in subdomains bounded by discontinuity lines (DLs), 
the latter defined from the Rankine-Hugoniot jump conditions.
The proposed approach allows to efficiently consider an arbitrary number of entropic shock waves, shock wave generation, as well as wave interactions. Some numerical experiments are presented to illustrate the strengths and properties of the algorithms.
\end{abstract}


\end{frontmatter}

%\tableofcontents

\section{Introduction}\label{sec:intro}
 In this paper we propose to develop a non-diffusive neural network (NDNN) solver for the accurate computation of shock waves in nonlinear (quasi-linear) hyperbolic conservation laws (HCLs).
 The objective is to track sharply entropic shocks and to circumvent a well-known issue when approximating HCL with Physics informed neural network (PINN) methods, more specifically when approximating shock waves with neural networks. 
 We consider the following Initial Value Problem (IVP):
 \begin{subequations}\label{e:HCL0}
 \begin{align}
 \mbox{\it for $f \in C^2$ convex}
 &
 \;\;\mbox{\it and 
 $u_0 \in BV(\Omega)$,
% $u_0 \in L^1(\R)\cap \textrm{BV}(\R)$,
 	 $\Omega\subseteq\R$, find $u$:}
\nonumber \\
\partial_t u + \partial_{x}f(u) 
& = 0, \;\;\; \mbox{\it in}\;\;\; Q:=\Omega\times (0,T],	\label{e:HCL0;PDE}\\
u(\cdot,0)&=u_0,\;\; \mbox{\it on}\;\; \Omega.		\label{e:HCL;bc}
\end{align}
 \end{subequations}
\textcolor{black}{In addition for $\Omega$ bounded, we impose boundary conditions for incoming characteristics, see for instance \cite{jmg,god1}.}
 We refer \cite{serre,lefloch,smoller} for the analysis of HCLs, and to \cite{god1,god2,despres} for their standard numerical approximation using finite volume$/$difference methods. 
 
\textcolor{black}{In this paper, we do not consider non-convex flux functions. The latter is known to generate non-classical shocks, see \cite{lefloch}. However ideas similar to those developed in this paper could be developed for entropic non-classical shocks, by combining piecewise smooth functions and Rankine-Hugoniot solvers (as in the convex case for first order ODE), and additional equations derived from the weak formulation of the PDE specific to non-convex flux \cite{laforest2008diminishing}.}

\subsection{Introductory remarks}
We recall that PINN algorithms allow for the computation of solutions to partial differential equations (PDE), and corresponding inverse problems, by using parameterized neural networks. More specifically, 
the solution $u$ is approximated by a parameter-dependent network $N$, and the $L^2$-norm (other norms can be used) of the residual of the equation applied to $N$  is minimized by standard (stochastic) gradient descent-type methods. 
One of the main strengths of this approach, which was originally developed in its most simple form by Lagaris \cite{lagaris}, is the use of automatic differentiation of explicit neural networks. There is no need for approximating differential operators, hence avoiding to a certain extent stability issues for evolution PDE. The computation of direct and inverse PDE problems with neural networks has become a very active research area from the practical, numerical as well mathematical points of view. 
Notice however that a full mathematical analysis of convergence, accuracy and stability is still far from complete. We refer to \cite{pinns,pinns2,pinns3} for details. Notice hereafter that other types of neural network-based algorithms have been developed, see \cite{otherNN,otherNN2,otherNN3}.

The approximation of shock waves using direct PINN \cite{pinns} algorithms can be inaccurate/very diffusive, or  even simply not convergent \cite{HypPINN}. 
\textcolor{black}{Among recent papers devoted to the numerical computation of HCL using neural networks, let us mention \cite{BLP} where is proposed a physics-informed attention-based neural network (PIANN) for non-convex fluxes generating non-classical shock waves \cite{lefloch}. As in our work here, the neural networks in \cite{BLP} are designed to include 
%\textcolor{blue}{from??} 
some knowledge of the structure of the solution. In \cite{PATEL2022110754} is proposed a least squares space-time control volume scheme, using space-time integral form with strong connection with finite volume methods.} 

By default, neural networks are constructed using smooth activation functions,  
which precludes the construction of weak discontinuous waves.
However, one has the freedom to choose also non-differentiable (ReLU, for instance) activation
functions. In this case, the direct application of PINN leads to differentiating a non-smooth function. 
 \textcolor{black}{We note that ReLU functions could be used as the activation function, as long as the training points do not coincide with the points where ReLU functions are not differentiable. However, the main issue in approximating shock waves, comes from that as weak solutions, shock waves are described mathematically from a weak formulation (leading to Rankine-Hugoniot jump conditions), hence independently of the regularity of the chosen activation function. In conclusion, we do not claim that it is impossible to use ReLU functions, but it should be done very carefully and a priori, on the weak formulation of the PDE.}
In addition, for complex solutions containing several waves (shocks, rarefactions) the convergence of the loss function is difficult to achieve. Notice however that the approximation of rarefaction can usually be accurately performed using deep and complex enough networks. 

In the present paper, the proposed approach is focused on the accurate approximation of shock waves. It allows for the generation of shock waves, as well as shock-shock or rarefaction-shock interactions.
The principle of the proposed methodology is: i) to use space-time dependent neural networks to solve HCLs in space-time subdomains bounded by curves of discontinuities; ii) use time-dependent neural networks to identify the discontinuity lines (DLs),
iii) solve HCL in space-time subdomains where the corresponding solutions are smooth and identify the DLs by 
minimizing a loss functional measuring the HCL residual in all subdomains and 
Rankine-Hugoniot's jump conditions on all the DLs.

More specifically, we decompose the  global space-time domain $Q$ in subdomains, 
which are delimited by the DLs defining shock waves   which initially are unknown. The solution in each subdomain, which are bounded by the DLs, and the DLs are then represented by dedicated neural networks.
The networks are trained by minimizing a loss functional, which measures the error of the PDE 
\eqref{e:HCL0;PDE} in each subdomain (for training the subdomain networks),
the error for approximating the initial condition \eqref{e:HCL;bc},
and	the error of Rankine-Hugoniot conditions (for training the networks dedicated to DLs).

The loss functional is nonlinear, and the associated minimization problem is not easy.
In the machine learning community, these problems are typically solved with a gradient descent method, or variants of it. In this paper we use a global gradient method.
We propose also a domain decomposition method (DDM) for the minimization problem, which 
allows the decoupling of the computation of the local solutions and DLs, hence allowing for an embarrassingly parallel computation, see \cite{jcp2022,NA2023}.

\textcolor{black}{Notice that in this paper we focus on one-dimensional HCLs.}
We present a proof of concept, as well as some analytical arguments justifying the application or extension to high dimensional problems. In future works, we plan to apply the derived methodology to high dimensional problems.

\subsection{Basics of neural networks}
In this subsection, we recall some basic concepts on neural networks. 
From the scientific machine learning point of view, neural networks are nothing but the composition of explicit parameterized linear functions built from discrete convolution products and of nonlinear (smooth or not) activation functions. 
More precisely, let $d$ be the dimension of the space,
	$l\in\mathbb N$, $n_i\in\mathbb N$, $i=0, 1,\ldots, l$, with $n_0=d$.
	For $\theta \in  \prod_{i=1}^l \mathbb R^{n_i\times (n_{i-1}+1)}$ we write
	\begin{eqnarray*}
		\theta&=&(\theta^i),\quad \theta^i=(w^i,b^i),\quad i=1,\ldots,l,\;\;\mbox{\it where}\\
		&&w^i(\cdot,\cdot)\in \mathbb R^{n_i\times n_{i-1}},\\
		&&b^i(\cdot)\in \mathbb R^{n_i}.
	\end{eqnarray*}
	A (fully connected) network in $\mathbb R^d$ with architecture $A=[n_0,n_1,\ldots,n_l]$ is a function of the form
	\begin{subequations}
	\begin{align}
	N_A&:(\theta;x)\in \prod_{i=1}^l \mathbb R^{n_i\times (n_{i-1}+1)}\times \mathbb R^d 
	\mapsto
	N_A(\theta;x)\in \mathbb R^{n_l},
	\\
	N_A(\theta;x)
	&=
	w^l\cdot\sigma(w^{l-1}\cdot(\cdots \sigma(w^2\cdot\sigma(w^1\cdot x + b^1)+b^2)\cdots)+b^{l-1})+b^l,
	\end{align}
	\end{subequations}
	\textcolor{black}{where  $\sigma:\mathbb R\mapsto\mathbb R$ is a given function, called activation function}, which acts on any vector or matrix component-wise.
It is clear from this exposition that the network $N_A$ is defined uniquely by 
the architecture $A$. 
In all the networks we will consider the architecture $A$ is given/fixed,
and we will omit the letter $A$ from $N_A$.

These network functions, which will be used to approximate solutions to HCLs, benefit from automatic differentiation with respect to $x$ and ${\theta}$ (parameters). This feature allows an evaluation of \eqref{e:HCL0} without error.
Naturally, the computed solutions are restricted to the function space spanned by the neural networks.

In the following, we will use networks in dimension $1+1$,
for approximating the solution of HCLs in space-time subdomains (one dimensional space).
We will denote these networks by latin bold capital letters, such as 
${\bf N}_i$, and by $\bs{\Theta}_{i}$ the associated parameters.
For such networks we will write ${\bf N}_i(\bs{\Theta}_{i};x,t)$, or whenever there is no ambiguity, 
simply ${\bf N}_i(x,t)$.

Similarly, we will use networks in dimension $0+1$,
for approximating the discontinuity lines.
We will denote these networks by latin bold lowercase letters, such as 
${\bf n}_i$, and by $\bs{\uptheta}_{_i}$ the associated parameters.
For such networks we will write ${\bf n}_i(\bs{\uptheta}_{i};t)$, or whenever there is no ambiguity, 
simply ${\bf n}_i(t)$.

\subsection{About the entropy of direct PINN methods}\label{subsec:SW5}
We finish this introduction by a short discussion on the entropy of PINN solvers for HCLs.  As direct PINN solvers can not directly capture shock waves (discontinuous weak solution), artificial diffusion is usually added to HCLs, $u_{\varepsilon}$ is searched as the solution to
\begin{eqnarray*}
	\partial_t u_{\varepsilon} + \partial_{x}f(u_{\varepsilon} )
	 & = & \varepsilon \partial_x(b(u_{\varepsilon})\partial_x u_{\varepsilon}) ,\\
	u_\varepsilon&=&u_0,\quad\mbox{\it with\;\;\; $u_0(x)\xrightarrow[x\to\pm L]{} u_{\pm}$},
\end{eqnarray*}
where $\varepsilon=o(1)$, $b>0$ and $\|b\|_{L^\infty}<\infty$.   Let us recall that entropic shock waves satisfy the Rankine-Hugoniot jump condition between the left and right states of a discontinuity $(u^{-},u^{+})$, i.e. $s(u^{-}-u^{+})=f(u^{-})-f(u^{+})$, where $s$ is the speed of the shock wave and $f'(u^{+})<s<f'(u^{-})$.  Unlike the case $\varepsilon=0$, the regularized equation has a unique smooth solution. Interestingly and at least in the scalar case, it is observed that direct PINN solvers compute ``entropic'' viscous shock profiles or rarefaction waves, but not ``non-entropic'' viscous shock profiles: if $f'(u_{-})>f'(u_{+})$ the PINN algorithm approximates viscous shock profiles $(x,t) \mapsto v\big((x-\sigma t)/\varepsilon\big)$, which theoretically converges in the distributional sense to an entropic shock,  and to a rarefaction wave for $f'(u_{+})>f'(u_{-})$. 

While we propose in Experiment 6 an illustration of this property of the PINN method,
% In this sense, we can consider direct PINN algorithm as an ``entropic solver'' although practically it only computes classical solutions. 
in this paper, we propose a totally different strategy which does not require the addition of any artificial viscosity, and computes entropic shock wave without diffusion.

\subsection{Organization of the paper}
This paper is organized as follows. Section \ref{sec:1d-solver} is devoted to the derivation of the basics of the non-diffusive neural network (NDNN) method. Different situations are analyzed including multiple shock waves, shock wave interaction, shock wave generation, and the extension to systems. In Section \ref{sec:ddm}, we discuss the efficient implementation of the derived algorithm. In Section \ref{sec:numerics}, several numerical experiments are proposed to illustrate the convergence and the accuracy of the proposed algorithms. We conclude in Section \ref{sec:conclusion}.

\section{Non-diffusive neural network solver for one dimensional scalar HCLs}\label{sec:1d-solver}
In this section, we derive a NDNN algorithm for solving solutions to HCL containing (an arbitrary number of) shock waves. The derivation is proposed in several steps:
\begin{itemize}
\item Case of one shock wave.
\item Case of an arbitrary number of shock waves.
\item Generation of shock waves.
\item Shock-shock interaction.
\item Extension to systems.
\end{itemize}
Generally speaking for $D$ shock waves, the basic approach will require $2D+1$ neural networks, more specifically, $D+1$ two-dimensional neural networks for approximating the local HCL smooth solutions, and $D$ one-dimensional networks for approximating the DLs. 
%In the following, we will approximate the equation on a bounded space-time domain $\Omega %\times (0,T)$.

\subsection{One shock wave}\label{subsec:SW1}
%In this subsection, we present the key idea for diffusion-free$/$discontinuity-tracking computation within PINN %algorithms. 
%The numerical solution is constructed as a piecewise smooth weak solution to HCL. 
We consider \eqref{e:HCL0} with boundary conditions when necessary (\textcolor{black}{incoming characteristics at the boundary}),
with  $\Omega=(a,b)$, and (without restriction) $0\in\Omega$.
We assume that the corresponding solution contains an entropic shock wave, with discontinuity line (DL) initially located at $x=0$, parameterized by $\gamma:t \mapsto \gamma(t)$, 
with $t\in[0,T]$ and $\gamma(0)=0$.
The DL $\gamma$ separates $\Omega$ in two subdomains $\Omega^-$, $\Omega^+$ (counted from left to right) and 
$Q$ in two time-dependent subdomains denoted 
$Q^-$ and $Q^+$ (counted from left to right).  We denote by $u^\pm$ the solution $u$ of \eqref{e:HCL0} in $Q^\pm$.
Then \eqref{e:HCL0} is written in the form of a system of two HCLs
\begin{subequations}\label{e:HCL,Qpmg}
\begin{align}
\partial_t u^{\pm} + \partial_x f(u^{\pm}) 
&=  0,\hspace{6mm}\textit{ in }\;  Q^{\pm}\, ,\\
u^{\pm}(\cdot,0) 
& =  u_{0},\quad\textit{ on }\; \Omega^\pm\, ,
\end{align}
\end{subequations}
which are coupled through the Rankine-Hugoniot (RH) condition along the DL for $t \in [0,T]$,
\begin{eqnarray}
\gamma'(t)\big[u^+(\gamma(t),t) - u^-(\gamma(t),t)\big]&=&f(u^+(\gamma(t),t)) - f(u^-(\gamma(t),t))\, ,
\label{e:RH-1d}
\end{eqnarray}
with the Lax shock condition reading as
\begin{eqnarray}
    f'(u^-(\gamma(t),t)) > \gamma'(t) >  f'(u^+(\gamma(t),t)) \, . 
\label{e:Lax-1d}
\end{eqnarray}
The proposed approach consists in approximating the DL and the solutions in each subdomain $Q^\pm$
with neural networks. 
We denote by ${\bf n}(t)$ the neural network approximating $\gamma$, with parameters $\bs{\theta}$ and ${\bf n}(0)=0$. We will refer still by ${\bf n}$ to DL given by the image of ${\bf n}$.
Like in the continuous case, ${\bf n}$ separates $Q$ in two domains, which are denoted again by $Q^\pm$.
We denote by ${\bf N}^\pm(x,t)$ the neural networks with parameters $\bs{\Theta}^\pm$  approximating $u^{\pm}$ in $Q^\pm$, and by $\partial_x{\bf N}^\pm$, resp. $\partial_t{\bf N}^\pm$, its $x$, resp. $t$ derivatives.

It is useful to consider the domains $Q^\pm$ as the image of the 
transformations ${\bf T}^\pm$, each of them defined in a fixed rectangle $R:=(0,1)\times(0,T)$, defined by
\begin{eqnarray}
  {\bf T}^-
&:&R \rightarrow Q^-,
\qquad
{\bf T}^-(x,t)=\left(({\bf n}(t) - a)x + a,t\right),		\label{e:R-gamma} 
\\
  {\bf T}^+
  &:&R\rightarrow Q^+,
  \qquad
  {\bf T}^+(x,t)=\left((b-{\bf n}(t))x + {\bf n}(t),t\right).		\label{e:R+gamma}  
\end{eqnarray}
Hence,
\begin{equation}
Q^+={\rm Im}({\bf T}^+),
\quad
Q^-={\rm Im}({\bf T}^-),
\quad
Q=Q^+\cup{\bf n}\cup Q^-,		\label{Q=Q+cupgcupQ-}
\end{equation}
where we have identified ${\bf n}$ with its image.
Equations \eqref{e:HCL,Qpmg} can be written in $Q^\pm$ in the form of a system of two HCLs
\begin{subequations}\label{e:HCL,Qpmn}
\begin{align}
(\partial_t {\bf N}^{\pm}(\cdot) 
+ 
\partial_x f({\bf N}^{\pm}(\cdot))) \circ {\bf T}^\pm(x,t)& = 0,\hspace{11mm}(x,t)\in R,
\\
%
%\partial_t {\bf N}^{\pm}({\bf T}^\pm(x,t)) 
%+ 
%\partial_x f({\bf N}^{\pm}({\bf T}^\pm(x,t))) &= 0,\hspace{11mm}(x,t)\in R\, ,\\
{\bf N}^{\pm}({\bf T}^\pm(x,0)) &=  u_{0}(x),\quad x\in(0,1).
\end{align}
\end{subequations}
The RH conditions \eqref{e:RH-1d} are expressed in terms of ${\bf N}^\pm(x,t)$  and ${\bf n}(t)$ are written as
\begin{eqnarray}\label{e:RH,g}
  \partial_t {\bf n}(t)
\big[{\bf N}^+({\bf n}(t),t) - {\bf N}^-({\bf n}(t),t)\big]
&=&
f({\bf N}^+({\bf n}(t)))- f({\bf N}^-({\bf n}(t))),\quad t\in[0,T].
\end{eqnarray}
In general,  ${\bf N}^\pm$ and ${\bf n}$ only approximate \eqref{e:HCL,Qpmn} and \eqref{e:RH,g} Denoting 
  $\bs{\Theta}
  =
  (\bs{\Theta}^-,\bs{\Theta}^+,\bs{\uptheta})$, 
  we consider the following minimization problem:
 \begin{equation}
  \mbox{\it find }\;
  \bs{\Theta}_*\in\bs{\Theta}_{ad}\;\; \mbox{\it  such that }\;\;
  	\bs{\cal L}(\bs{\Theta}_*)=\min\{\bs{\cal L}(\bs{\Theta}),\;\; 
  	\bs{\Theta}\in\bs{\Theta}_{ad}
  \},
  \label{e:min(L)}
\end{equation}
where
\begin{eqnarray}
\bs{\cal L}(\bs{\Theta}) 
 & = & 
    \lambda 
    \Big(
     \big\|(\partial_t {\bf N}^-(\cdot) + 
      \partial_x f({\bf N}^-(\cdot)))\circ{\bf T}^- \|^2_{L^2(R)} +
     \nonumber\\
&&
\hspace*{5mm}
    \|(\partial_t {\bf N}^+(\cdot) + 
      \partial_x f({\bf N}^+(\cdot)))\circ{\bf T}^+ \|^2_{L^2(R)}  
    \Big)
    \nonumber\\
	& & + 
	\mu
	\big\|\partial_t {\bf n}(t)
	\big[{\bf N}^+({\bf n}(t),t) - {\bf N}^-({\bf n}(t),t)\big] - \nonumber\\
	&&
	\hspace*{19mm}
	\big[f({\bf N}^+({\bf n}(t),t))- f({\bf N}^-({\bf n}(t),t))\big]\|_{L^2(0,T)}^2
	\nonumber\\
    & & +
    \kappa
    \Big(
    \big\|{\bf N}^-({\bf T}^-(x,0)) - u_0({\bf T}^-(x,0))\|^2_{L^2(0,1)}  + 
    \nonumber\\
   &&
	\hspace*{9mm}
     \big\|{\bf N}^+({\bf T}^+(\cdot,0)) - u_0({\bf T}^+(\cdot,0))\|^2_{L^2(0,1)} 
    \Big),
      \label{e:L,1sw}
    \end{eqnarray}
for some positive parameters $\lambda,\mu$ and $\kappa$, and $\bs{\Theta}_{ad}$ is the set of admissible weights.

Let 
$\bs{\Theta}_*=
  (\bs{\Theta}^-_*,\bs{\Theta}^+_*,\bs{\uptheta}_*)$
be the solution of \eqref{e:min(L)}.
Then $\gamma$ is approximated by the network with parameters $\bs{\uptheta}_*$.
The network ${\bf n}_*$ divides $Q$ in two  domains $Q^-_*$ and $Q^+_*$. 
The solution $u^\pm$ in each of these subdomains is approximated by the networks
${\bf N}^\pm_*$ with parameters $\bs{\Theta}^\pm_*$.

\subsection{Arbitrary number of shock waves}\label{subsec:SW2}
We consider again \eqref{e:HCL0} with  $\Omega=(a,b)$, and boundary conditions  when necessary.
We assume that the solution is initially constituted by:
i) $D$ entropic shock waves, 
ii) an arbitrary number of rarefaction waves, and that 
iii) there is no shock generation for $t\in [0,T]$. 
We assume that the $D$ discontinuity lines are initially located at $x_i$, $i=1,\ldots,D$,
which motivates the decomposition of  the domain $\Omega$ in $D+1$ 
subdomains, $\Omega = \cup_{i=1}^{D+1}\Omega_i\cup_{i=1}^D\{x_i\}$, 
$\Omega_i=(x_{i-1},x_i)$, $x_0=a$, $x_{D+1}=b$.

For $i \in \{1,\ldots,D\}$, we denote by $\gamma_i:t\rightarrow \gamma_i(t)$ the DLs such that $\gamma_i(0)=x_i$ and by $\gamma'_i(t)$ the corresponding shock velocity. 
The DLs divide $Q$ in time-dependent domains $Q_i$. Namely, $Q_i$ is the subdomain of $Q$ bounded by the DLs $\gamma_{i-1}$ and $\gamma_i$, $i=1,\ldots,D+1$. 
We note that it is practical to denote $Q_i^-=Q_{i-1}$, $Q_i^+=Q_i$, $i=1,\ldots,D$, and 
$u_i^\pm$ the solution of \eqref{e:HCL0} in $Q_i^\pm$.
Then \eqref{e:HCL0}  is written equivalently as a system of $(D+1)$ HCLs 
\begin{subequations}\label{e:HCL,Qi;D}
	\begin{align}
	\partial_t u_i + \partial_x f(u_i) 
	&=  0,\hspace{11mm}\textit{ in }\;  Q_i\, ,\\
	u_i(\cdot,0) 
	& =  u_{0}(\cdot),\quad\textit{ on }\; \Omega_i\, ,\quad i=1,\ldots,D+1,
	\end{align}
\end{subequations}
which are coupled through the RH conditions along the DLs for $t \in [0,T]$,
\begin{eqnarray}
\gamma'_i(t)\big[u_i^+(\gamma_i(t),t) - u_i^-(\gamma_i(t),t)\big]
&=&
f(u_i^+(\gamma_i(t),t)) - f(u_i^-(\gamma_i(t),t))\, ,
\label{e:RH-1d;D}
\end{eqnarray}
with Lax entropy condition satisfied
\begin{eqnarray}
f'(u_i^-(\gamma_i(t),t))>\gamma'_i(t)>f'(u_i^+(\gamma_i(t),t)),\quad i=1,\ldots,D.
\label{e:Lax-1d;D}
\end{eqnarray}

Similar to the previous case, the proposed approach consists in approximating the solutions in each subdomain $Q_i$ and the DL with neural networks. 
We denote by ${\bf n}_i(t)$ the neural network approximating $\gamma_i$, with parameters $\bs{\uptheta}_i$
and ${\bf n}_i(0)=x_i$.
Like the DLs $\gamma_i$, 
${\bf n}_i$ separate $Q$ in $D+1$ domains, which are denoted again by $Q_i$.
We denote by ${\bf N}_i(x,t)$ the neural networks with parameters 
$\bs{\Theta}_i$ approximating $u$ in $Q_i$ . 

Like in the case of one shock wave in Section 
\ref{subsec:SW1}, it is useful to consider $Q_i$ as the image of the transformations ${\bf T}_i$
defined as follows. For $i=1,\ldots,D+1$ define
\begin{eqnarray}
{\bf T}_i
&:&R\rightarrow Q_i,
\qquad
{\bf T}_i(x,t)=\left(({\bf n}_i(t)-{\bf n}_{i-1}(t))x + {\bf n}_{i-1}(t),t\right),
			\label{e:Ri}  
\end{eqnarray}
with ${\bf n}_0=a$ and ${\bf n}_{D+1}=b$.
Hence,
\begin{equation}
Q_i={\rm Im}({\bf T}_i),
\quad
Q=\cup_{i=1}^{D+1} Q_i\cup_{i=1}^D{\bf n}_i.					\label{e:Q=cupQi}
\end{equation}
Like in Section \ref{subsec:SW1}, it is convenient to denote $Q_i^-=Q_i$ and $Q_i^+=Q_{i+1}$, $i=1,\ldots,D$.  Then \eqref{e:HCL,Qi;D} and \eqref{e:RH-1d;D} are written in terms of ${\bf N}_i$ and ${\bf n}_i$ as follows
\begin{subequations}\label{e:HCL,Qi,n}
	\begin{align}
	(\partial_t {\bf N}_i(\cdot) 
	+ 
	\partial_x f({\bf N}_i(\cdot))\circ{\bf T}_i(x,t) &= 0,\hspace{23mm}\textit{ in }\;  R \, ,\\
	{\bf N}_i({\bf T}_i(x,0)) &=  u_{0}({\bf T}_i(x,0)),\quad\textit{ on }\; (0,1) \, ,
	\end{align}
\end{subequations}
and the Rankine-Hugoniot conditions
\begin{eqnarray}\label{e:RH,i}
\hspace*{-7mm}
\partial_t {\bf n}_i(t)
\big[{\bf N}_i^+({\bf n}_i(t),t) - {\bf N}_i^-({\bf n}_i(t),t)\big]
&=&
f({\bf N}_i^+({\bf n}_i(t),t))- f({\bf N}_i^-({\bf n}_i(t),t)),\, t\in[0,T].
\end{eqnarray}
In general, ${\bf N}_i^\pm$ and ${\bf n}_i$ do not exactly solve \eqref{e:HCL,Qi,n} and \eqref{e:RH,i}.
Denoting 
$\bs{\Theta}
=
\prod_{i=1}^{D+1}\bs{\Theta}_{i}\times\prod_{i=1}^D\bs{\uptheta}_{i}$, 
the optimized networks (or equivalently parameters) are obtained by solving the problem:
\begin{equation}
\mbox{\it find }\;
\bs{\Theta}^*\in\bs{\Theta}_{ad}\;\; \mbox{\it  such that }\;\;
\bs{\cal L}(\bs{\Theta}^*)
=
\min\{\bs{\cal L}(\bs{\Theta}),\;\; 
\bs{\Theta}\in\bs{\Theta}_{ad}
\},
\label{e:min(L),D}
\end{equation}
where
\begin{eqnarray}
\bs{\mathcal L}(\bs{\Theta}) 
& = & 
\lambda 
\sum_{i=1}^{D+1}
\big\|
(\partial_t {\bf N}_i(\cdot) 
+
\partial_x f({\bf N}_i(\cdot))
)\circ{\bf T}_i\|^2_{L^2(R)} 
\nonumber\\
&&+
\mu
\sum_{i=1}^D
\big\|\partial_t {\bf n}_i(\cdot)
\big[{\bf N}_i^+({\bf n}_i(\cdot),\cdot) - {\bf N}_i^-({\bf n}_i(\cdot),\cdot)\big] - \nonumber\\
&&
\hspace*{27mm}
\big[f({\bf N}_i^+({\bf n}_i(\cdot),\cdot))- f({\bf N}_i^-({\bf n}_i(\cdot),\cdot))\big]\|_{L^2(0,T)}^2
\nonumber\\
& & + 
\kappa
\sum_{i=1}^{D+1}
\|{\bf N}_i({\bf T}_i(\cdot,0)) - u_0({\bf T}_i(\cdot,0))\|^2_{L^2(0,1)},
\label{e:L,Dsw}
\end{eqnarray}
for some positive parameters $\lambda,\mu$ and $\kappa$,
and where $\bs{\Theta}_{ad}$ is the set of admissible weights.

Like in Section \ref{subsec:SW1}, if 
$\bs{\Theta}_*=\prod_{i=1}^{D+1}{\Theta}^*_i\times\prod_{i=1}^D\bs{\uptheta}^*_i$
be the solution of \eqref{e:min(L)} then
the DLs $\gamma_i$ are approximated by the networks ${\bf n}^*_i$ with parameters $\bs{\uptheta}^*_i$. These  networkss divide $Q$ in $D+1$ subdomains denoted by $Q^*_i$.
The solution in each of these subdomains is approximated by the networks 
${\bf N}^*_i$ with parameters $\bs{\Theta}^*_i$.

This approach involves $2D+1$ neural networks. 
Hence the convergence of the optimization algorithm may be hard to achieve for large $D$. 
However, the solutions $u_i$ being smooth and DLs being dimensional functions,
the networks ${\bf N}_i$ and ${\bf n}_i$ do not need to be deep.

\subsection{Shock wave generation}\label{subsec:SW3}
So far we have not discussed the generation of shock waves within a given domain. 
Interestingly the method developed above for pre-existing shock waves can be applied directly for the generation of shock waves in a given subdomain. 
In order to explain the principle of the approach, we simply consider one space-time domain 
$\Omega \times [0,T]$. We initially assume that:
i) $u_0$ is a smooth function, and 
ii) at time $t^* \in (0,T)$ a shock is generated in $x^{*} \in \Omega$, and the corresponding 
DL is defined by $\gamma: t \mapsto \gamma(t)$ for $t\in [t^*,T]$. 

\textcolor{black}{Although $t^*$ could be analytically estimated from $-1/\min_x f'(u_0(x))$, our algorithm does not require {\it a priori} the knowledge of $t^*$ and $x^*$. The only required information is the fact that a shock will be generated which can again be deduced by the study of the variations of $f'(u_0)$.} 

Unlike the framework in Sections 2.1 and 2.2, here we cannot initially specify the position of the DL. However proceed here using a similar approach.
Let $x_0\in\Omega$ be such that $\gamma$ can be extended as a smooth curve in $[0,T]$, still denoted by $\gamma$, 
with $\gamma(0)=x_0$. 
Without loss of generality we may assume $x_0=0$.
Then we proceed exactly as in Section \ref{subsec:SW1}, with one network ${\bf n}$ representing the DL, 
and ${\bf N}^\pm$ the solutions on each side of ${\bf n}$. 
We define $t^*=\max\{t\in(0,T],\;\; {\bf N}^+({\bf n}(t),t)={\bf N}^-({\bf n}(t),t)\}$.
Note that for $t<t^*$ we have ${\bf N}^+({\bf n}(t),t)={\bf N}^-({\bf n}(t),t)$, 
and the curve $\{{\bf n}(t),\; t\in[0,t^*)\}$ does not have any significant meaning.

\subsection{Shock wave interaction}\label{subsec:SW4}
As described above, for $D$ pre-existing shock waves we decompose the global domain in $D+1$ subdomains. So far we have not considered shock wave interactions leading to a new shock wave, which reduces by one the total number of shock waves for each shock interaction. 

Assume that two interacting shock waves with DLs $\gamma_{i-1}(t)$, $\gamma_i(t)$ are managed through subdomains $Q_{i-1}$, $Q_{i}$, $Q_{i+1}$,
and intersect at $t=t^*$ where
$\gamma_{i-1}(t^*)=\gamma_i(t^*)=x^*$, for a certain $x^*$.
If we set $\Omega_i^t=\{{\bf T}_i(x,t),\; x\in(0,1)\}$, see Section \ref{subsec:SW2} for the notations,
it is expected that the domain $\Omega_i^t$ becomes empty at $t^*$.
In this case, we proceed as follows:
\begin{itemize}
\item  Solve \eqref{e:HCL0} in $[0,T^*]$, where $T^*$ is an estimated time of interaction such that $T^*>t^*$ and close to $t^*$. 
\item 
Deduce precisely $t^*$, where $\gamma_{i-1}(t^*)=\gamma_i(t^*)$.
\item 
Re-decompose the global domain $\Omega$ based on the fact that 
at time $t^*$, two chock waves interact at $x^*=\gamma_{i-1}(t^*)=\gamma_i(t^*)$.
\item 
Compute the solution for $t>t^*$.
\end{itemize}
The evaluation of $T^*$ can be performed by linearizing the system or by restart; once $t^*$ is accurately computed the global domain can be re-decomposed.

\subsection{Non-diffusive neural network solver for one dimensional systems of CLs}\label{sec:nd-solver}
In this subsection, we extend the above ideas to hyperbolic systems of conservation laws. Let us denote  $f=(f_1,\ldots,f_m)\in C^2(\R^m;\R^m)$,  such that $A(u)=\Big[\partial_{u_j}f_i(u)\Big]$, $u=(u_1,\ldots,u_m)\in\R^m$, is strictly hyperbolic, and  
%$u^0=(u^0_1,\ldots,u^0_m) \in L^1(\R;\R^m)\cap \textrm{BV}(\R;\R^m)$, 
$u^0=(u^0_1,\ldots,u^0_m) \in \textrm{BV}(\R;\R^m)$, 
$\Omega\subseteq\R$.  We look for a solution $u=(u_1,\ldots,u_m)$ to the following initial value problem
\begin{subequations}\label{pde+bc;d}
\begin{align}
\partial_t u + \partial_{x}f(u) & = 0,\;\;\;\;\;\;\;\,
\hbox{\it in } Q:=\Omega\times (0,T),\quad \Omega=(a,b),		\label{pde,Rx(0,T);d}\\
u(\cdot,0)&=u^0(\cdot),\;\, \hbox{\it on } \Omega.						\label{bc,R;d}
\end{align}
\end{subequations}
We refer \cite{serre,lefloch,smoller} for the analysis of hyperbolic systems of conservation laws, and to \cite{god1,god2} for their standard numerical approximation using finite volume$/$difference methods. The extension of the method developed above is in principle straightforward. In the following we consider piecewise smooth solutions to \eqref{pde+bc;d}.

If $\gamma(t)$ describes a DL and $u^-$, resp. $u^+$, is the solution on the left, resp. right, of $\gamma$,
as in \eqref{e:RH-1d} we have
\begin{eqnarray}
\gamma'(t)
\big[u^+(\gamma(t),t)) - u^-(\gamma(t),t)\big]
&=&
f(u^+(\gamma(t),t)) - f(u^-(\gamma(t),t)).	\label{e:RH-nd}
\end{eqnarray}
Moreover the Lax shock conditions read as follow for $k\in\{1,\ldots,m\}$: if the $k$th characteristic 
field is genuinely nonlinear, then
\begin{subequations}	\label{e:Lax-nd}
\begin{align}
	\lambda_{k}(u^+(\gamma(t),t)) &< \gamma'(t) < \lambda_{k+1}(u^+(\gamma(t),t)),\\
	\lambda_{k-1}(u^-(\gamma(t),t)) &< \gamma'(t) < \lambda_{k}(u^-(\gamma(t),t)),
\end{align}
\end{subequations}
and if it is linearly degenerate 
\begin{equation}	\label{e:Lax-nd2}
	\lambda_{k}(u^-(\gamma(t),t)) = \gamma'(t) = \lambda_{k}(u^+(\gamma(t),t)),
\end{equation}
where $\lambda_1(u)<\cdots<\lambda_m(u)$ are the eigenvalues of $A(u)$.

Unlike the scalar case, Riemann's problems for systems require an initial decomposition in up to $m$ (shock, contact discontinuity, rarefaction) waves.
Hence, if a DL emanates from $x_i\in\Omega$, $i=1,\ldots,D$, in fact there are up to $m$ DLs, which will be assumed in the presentation hereafter. We will denote these DLs by $\gamma_{i,j}$, $j=1,\ldots,m_i$, where $\gamma_{i,j}'(0)$ are ordered increasing in $j$ and $m_i\leq m$. 
We set $\gamma_{0,1}=a$, $m_0=1$, $\gamma_{D+1,1}=b$, $m_{D+1}=1$, and 
$\gamma_{i,0}=\gamma_{i-1,m_{i-1}}$, $i=1,\ldots,D+1$.
Then $\gamma_{i,j}$, $i=0,\ldots,D+1$, $j=1,\ldots,m_i$, separate $Q$ in  subdomains denoted 
$Q_{i,j}$,   $i=1,\ldots,D+1$, $j=1,\ldots,m_i$, where
$Q_{i,j}$ is bounded by $\gamma_{i,j}$ and $\gamma_{i,j-1}$. 
We also set $Q_{i,m_i+1}=Q_{i+1,1}$, $i=1,\ldots,D$.

If $u_{i,j}$ is the solution of \eqref{pde+bc;d} in $Q_{i,j}$, then \eqref{pde+bc;d}  can be equivalently written as 
\begin{subequations}\label{pde+bc;d;i,0}
	\begin{align}
	\partial_t u_{i,1} + \partial_{x}f(u_{i,1}) & = 0,\;\;\;\;\;\;\; 
	\hbox{\it in } Q_{i,1},		\label{pde,Rx(0,T);d;i,0}\\
	u_{i,1}(\cdot,0)&=u^0(\cdot),\;\, \hbox{\it on } \Omega_i:=(x_{i-1},x_i),\quad
	i=1,\ldots,D+1,						\label{bc,R;d,i,0}
	\end{align}
\end{subequations}
and
\begin{subequations}\label{pde+bc;d;i,j}
	\begin{align}
	\partial_t u_{i,j} + \partial_{x}f(u_{i,j}) & = 0,\;\;\;\;\;\;\, 
	\hbox{\it in } Q_{i,j},\quad i=1,\ldots,D,\;\; j=2,m_i.		\label{pde,Rx(0,T);i,j}
	%\\
%	u_{i,j}(\cdot,0)&=u^0(\cdot),\;\, \hbox{\it on } \Omega_{i},\quad
%	i=1,\ldots,D,\;\; j=2,m_i,						\label{bc,R;d;i,j}
	\end{align}
\end{subequations}
On $\gamma_{i,j}$, $i=1,\ldots,D$, $j=1,\ldots,m_i$, the Rankine-Hugoniot condition  \eqref{e:RH-nd} 
is written as
\begin{eqnarray}
\gamma'_{i,j}(t)
\big[u_{i,j}^+(\gamma_{i,j}(t),t) - u_{i,j}^-(\gamma_{i,j}(t),t)\big]
&=&
f(u_{i,j}^+(\gamma_{i,j}(t),t)) - f(u_{i,j}^-(\gamma_{i,j}(t),t)),
\label{e:RH-nd;i,j}
\end{eqnarray}
where $u_{i,j}^-=u_{i,j}$, and $u_{i,j}^+=u_{i,j+1}$.

The approximate solution and approximate DLs will be searched in the form of neural networks. Each DL $\gamma_{i,j}$ is approximated by a scalar network ${\bf n}_{i,j}$ with parameters 
$\bs{\uptheta}_{{i,j}}$ and ${\bf n}_{i,j}(0)=x_i$.
Also we set 
${\bf n}_{0,1}=a$, ${\bf n}_{D+1,1}=b$,  and
${\bf n}_{i,0}={\bf n}_{i-1,m_{i-1}}$, $i=1,\ldots,D+1$.

The DLs ${\bf n}_{i,j}$ divide $Q$ in subdomains, which we denote again by $Q_{i,j}$, where
$Q_{i,j}$ is bounded by $\gamma_{i,j}$ and $\gamma_{i,j-1}$, $i=1,\ldots,D+1$, $j=1,\ldots,m_i+1$. 
It is convenient to set $Q_{i,m_i+1}=Q_{i+1,1}$, $i=1,\ldots,D$.

We can write Equations \eqref{pde+bc;d;i,0}, \eqref{pde+bc;d;i,j} and \eqref{e:RH-nd;i,j} 
in terms of ${\bf N}_{i,j}$ and ${\bf n}_{i,j}$ as follows.
First, like in Sections \ref{subsec:SW1} and \ref{subsec:SW2} we consider $Q_{i,j}$ as the image of the transformations ${\bf T}_{i,j}$ defined as follows. 
If $R=(0,1)\times(0,T)$ (a rectangle), $C=\{(x,t),\;\; |x|\leq t\leq T\}$ (a cone at the origin),
for $i=1,\ldots,D+1$, we define
\begin{eqnarray}
{\bf T}_{i,1}
&:&
R\rightarrow Q_{i,1},\quad \;\;i=1,\ldots,D+1,
\nonumber\\
{\bf T}_{i,1}(x,t)
&=&
\left(({\bf n}_{i,1}(t)-{\bf n}_{i,0}(t))x 
+ 
{\bf n}_{i,0}(t),t\right),		\label{e:Ri0}  
\\
{\bf T}_{i,j}
&:&
C\rightarrow Q_{i,j},\quad \;\;i=1,\ldots, D,\; j=2,\ldots,m_i,
\nonumber\\
{\bf T}_{i,j}(x,t)
&=&
\left(
\frac{t-x}{2t}({\bf n}_{i,j-1}(t)-{\bf n}_{i,j}(t)) 
+ 
{\bf n}_{i,j}(t),t\right) .		\label{e:Rij}  
\end{eqnarray}
Hence,
\begin{equation}
Q_{i,j}={\rm Im}({\bf T}_{i,j}),
\quad
Q={\bigcup_{i=1}^{D+1}}\bigcup_{j=1}^{m_i} Q_{i,j}\bigcup_{i=1}^D\bigcup_{j=1}^{m_i}{\bf n}_{i,j}.					\label{e:Q=cupQij}
\end{equation}
Then in terms of ${\bf N}_{i,j}$ and ${\bf n}_{i,j}$, equations \eqref{pde+bc;d;i,0} are written as
\begin{subequations}\label{e:HCL,Qi0;N}
\begin{align}
(\partial_t {\bf N}_{i,1}(\cdot) 
+ 
\partial_x f({\bf N}_{i,1}(\cdot))\circ{\bf T}_{i,1}(x,t) &= 0,\hspace{26mm}\textit{ in }\;  R\, ,\\
{\bf N}_{i,1}({\bf T}_{i,1}(x,0)) &=  u_{0}({\bf T}_{i,1}(x,0)),\hspace{5mm}\textit{ on }\; (0,1),\;\;
i=1,\ldots,D+1,
\end{align}
\end{subequations}
the equations \eqref{pde+bc;d;i,j} are written as (for $i=1,\ldots,D$, $j=2,\ldots,m_i$)
\begin{equation}
	(\partial_t {\bf N}_{i,j}(\cdot) 
	+ 
	\partial_x f({\bf N}_{i,j}(\cdot)))\circ{\bf T}_{i,j}
	 = 0,\hspace{6mm}\textit{ in }\;  C,
	\label{e:HCL,Qij;N}
\end{equation}
and the equations \eqref{e:RH-nd;i,j} as (for $i=1,\ldots,D$, $j=1,\ldots,m_i$)
\begin{equation}
\partial_t {\bf n}_{i,j}(t)
\big[{\bf N}_{i,j}^+({\bf n}_{i,j}(t),t) - {\bf N}_{i,j}^-({\bf n}_{i,j}(t),t)\big]
=
f({\bf N}_{i,j}^+({\bf n}_{i,j}(t),t))- f({\bf N}_{i,j}^-({\bf n}_{i,j}(t),t)),	\label{e:RH,N}
\end{equation}
where ${\bf N}_{i,j}^-={\bf N}_{i,j}$ and ${\bf N}_{i,j}^+={\bf N}_{i,j+1}$.

In general, ${\bf N}_{i,j}$ and ${\bf n}_{i,j}$ do not solve 
\eqref{e:HCL,Qi0;N}, \eqref{e:HCL,Qij;N} and \eqref{e:RH,N} exactly, but only approximately.
Denoting 
$\bs{\Theta}
=
\prod_{i=1}^{D+1}\prod_{j=1}^{m_i}
\bs{\Theta}_{{i,j}}
\times
\prod_{i=1}^D\prod_{j=1}^{m_i}\bs{\uptheta}_{{i,j}}$, 
the optimized parameters $\bs{\theta}^*$ are obtained by solving the problem:
\begin{equation}
\mbox{\it find }\;
\bs{\Theta}^*\in\bs{\Theta}\;\; \mbox{\it  such that }\;\;
\bs{\cal L}(\bs{\Theta}^*)=\min\{\bs{\cal L}(\bs{\Theta}),\;\; 
\bs{\Theta}\in\bs{\Theta}_{ad}
\},
\label{e:min(L),D;m}
\end{equation}
where
\begin{eqnarray}
\bs{\cal L}(\bs{\Theta}) 
&= & 
\lambda 
\Big(\sum_{i=1}^{D+1}
	\|(\partial_t {\bf N}_{i,1}(\cdot) 
	+ 
	\partial_x f({\bf N}_{i,1}(\cdot))\circ{\bf T}_{i,1} \|^2_{L^2(R)}  
	+
	\nonumber\\
& & \hspace*{6mm}
    \sum_{i=1}^{D}\sum_{j=2}^{m_i}
	 \|(\partial_t {\bf N}_{i,j}(\cdot) + \partial_x f({\bf N}_{i,j}(\cdot))\circ{\bf T}_{i,j} \|^2_{L^2(C)}  
\Big)
 \nonumber\\
 & & + 
	\mu
	\sum_{i=1}^{D}
	\sum_{j=1}^{m_i}
	\|\partial_t {\bf n}_{i,j}(t)
	\big[{\bf N}_{i,j}^+({\bf n}_i(t),t) - {\bf N}_{i,j}^-({\bf n}_{i,j}(t),t)\big]
	-
	\nonumber\\
&&	\hspace*{34mm}
   \big[f({\bf N}_{i,j}^+({\bf n}_{i,j}(t),t))- f({\bf N}_{i,j}^-({\bf n}_{i,j}(t),t))\big]\|_{L^2(0,T)}^2
	\nonumber\\
& & +
\kappa
\sum_{i=1}^{D+1}
	\|{\bf N}_{i,1}({\bf T}_{i,1}(\cdot,0)) - u_0({\bf T}_{i,1}(\cdot,0))\|^2_{L^2(0,1)},
	\label{e:L,nsw2}
\end{eqnarray}
for some positive subdomain-independent parameters $\lambda,\mu$ and $\kappa$. 
\textcolor{black}{Notice that when the considered IVP problems involve cones$/$subdomains of ``very different'' sizes, the choice of the hyper-parameters may be taken subdomain-dependent. This question was not investigated in this paper.} If $m_i=1$ for all $i$, then the problem \eqref{e:min(L)}, with $\bs{\mathcal L}$ given by \eqref{e:L,nsw2}, is no different from the scalar case, except that \eqref{e:L,nsw2}  involves the evaluation of norms for vector valued functions. 

If $\bs{\Theta}^*$ is the solution of the problem \eqref{e:min(L),D;m}, then the approximates of $\gamma_{i,j}$ and of the solution
$u$ in each $Q_{i,j}$ is computed similarly as in Section \ref{subsec:SW2}.

\subsection{Efficient initial wave decomposition}\label{subsec:IWD}
In this section we present an efficient initial wave decomposition for a Rieman problem which can be used to solve the problem \eqref{pde+bc;d}.The idea is to solve a Riemann problem for  $t$ small, and once the states are identified we use NDNN method.  Here we propose an efficient neural network based approach for the initial wave decomposition, which can easily be combined with the NDNN method.  We first recall some general principles on the existence of at most $m$ curves connecting two constant states $u_L$ and $u_R$ in the phase space. We then propose a neural network approach for constructing the rarefaction and shock curves as well as their intersection.

\subsubsection{Initial wave decomposition for arbitrary $m$}
%\textcolor{black}
{We here recall some fundamental results on the initial wave decomposition in Riemann problems with $m\geq 1$, when $u_R$ is close enough to $u_L$.
Hereafter, we consider a $m$ equations system on a bounded domain $\Omega$ (and Dirichlet boundary conditions)
\begin{eqnarray}\label{ibvp2bis}
  \left.
  \begin{array}{rcl}
    \partial_t u + \partial_{x}f(u) & = & 0, \,\,\,\,\, \textrm{ on } \Omega\times [0,T], \\
    u(\cdot,0) & = & u^0, \,\,\, \textrm{ on } \Omega \, ,
  \end{array}
  \right.
\end{eqnarray}
such that for any $x \in \Omega$
\begin{eqnarray}\label{ibvp2bis2}
  u^0(x)=\left\{
  \begin{array}{ll}
u_L, & x<0 \, , \\
u_R, & 0<x \, .
    \end{array}
  \right.
\end{eqnarray}}
%\textcolor{black}
{ We assume that $|u_L-u_R|\ll1$, and that for all $k \in \{1,\ldots,m\}$ the $k$th characteristic field is either genuinely nonlinear or linearly degenerate. If $u_R$ is  in a neighborhood of $u_L$, then the Riemann problem   \eqref{ibvp2bis} has a unique weak solution constituted by at most $m+1$ states $u_L=u^*_0,u^*_{1},\cdots,u^*_{(\mu-1)},u^*_\mu=u_R$, $\mu\leq m$  separated by rarefaction waves, shock waves or contact discontinuities. Moreover the intermediate states are located at the intersection of simple curves in the phase space.} 
% \textcolor{black}
 {This can be summarized by the existence of a smooth mapping $\chi$ from $\R^m$ to $\R^m$ in a neighborhood of $0\in \R^m$ (see \cite{god1}) such that 
 %(J'AI PRECISE UN PEU PLUS... MAIS PAS TROP CAR C'EST DETAILLE DANS LE LIVRE DE GOD\&RAVIART)
\begin{eqnarray*}
\chi(\varepsilon) = \chi_\mu(\varepsilon_\mu;\chi_{\mu-1}(\varepsilon_{\mu-1};,\cdots;\chi_1(\varepsilon_1;u_L)\cdots))\, ,
  \end{eqnarray*}
where the mapping $\chi_k$ defines a $k$-wave, $\varepsilon=(\varepsilon_1,\cdots,\varepsilon_\mu) \in \R^m$, with $\chi(0)=u_L$ and  $\chi(\varepsilon)=u_R$. The intermediate states $u^*_1,\cdots,u^*_{\mu-1}$ are such that $u^*_{k+1}=\chi_k(\varepsilon_{k+1};u^*_k)$ for $k=0,\cdots,\mu-1$. We again refer to \cite{god1} for details. }

\textcolor{black}{Such a decomposition could be used to approximate the solution of \eqref{pde+bc;d} for $t$ small. Once the intermediate states connected by rarefaction$/$shock waves and contact discontinuity are identified, we can then apply the NDNN method derived in this paper. In the following,  we denote by $(\lambda_k(u), r_k(u))$, $k=1,\ldots,\mu$, the eigenpairs of  $\nabla_uf(u)$, i.e $\nabla_u f(u)\cdot r_k(u)=\lambda_k(u)r_k(u)$. For simplicity, we will omit the dependence of $(\lambda_k(u),r_k(u))$, and $\nabla f(u)$ in $u$. Below, we implement explicitly this idea for $m=2,3$.}
\subsubsection{Initial wave decomposition for 2 equation systems $(m=2,\, \mu=2)$}
%In this paragraph, the presentation is proposed for one-dimensional shock waves and $m=2$ equations. 
We search for $u=(u_1,u_2): \Omega\times (0,T]$ to $\R^2$ solution to \eqref{ibvp2bis}.
%\begin{subequations}\label{initD}
% \begin{align}
%\partial_t u_i + \partial_{x}f(u)  & = 0, \;\;\; \mbox{\it in}\;\;\; Q:=\Omega\%times (0,T],\\
%u(\cdot,0)&=u_0,\;\; \mbox{\it on}\;\; \Omega.		
%\end{align}
%\end{subequations}
%where $f=(f_1,f_2)$ is a smooth vector-valued convex function, and $u_0: \R\rightarrow \R^2$ is defined as follows (Riemann problem):
%
\\
\vspace{1mm}
\noindent
{\bf Generation of 2 shock waves.} We here consider the generation of $2$ shock waves from $u_0$. We detail the computation of the intermediate state via the construction and intersection of the shock curves (Hugoniot loci) in the phase-plane ($u_1,u_2$). Let us recall that  $k$-shock curves are defined as the following integral curves with $k=1,2$, see \cite{god1,god2}
\begin{eqnarray*}
  s'_k(\xi) & = & r_k(s_k(\xi)), \, \textrm{ with } \xi > 0 \, ,\\
  s_k(0) & = & s_k^0 \, ,
  \end{eqnarray*}
where $s^0_1=u_L$ and $s^0_2=u_R$. The searched intermediate state denoted here $u_1^*$,
  respectively defines a $1$-shock $(u_L,u_1^*)$ and a $2$-shock $(u_1^*,u_R)$ and which is such that for some 
$\xi^*>0$, $u_1^*=s_1(\xi^*)=s_2(\xi^*)$. Moreover, the corresponding $k$-shock speeds  are 
$\sigma_k = \lambda_k(s_{k}({\xi}^*))$ with $k=1,2$.
\\
\\
Following is the neural network strategy,
\begin{enumerate}
\item we optimize 2 vector-valued neural networks $\nu_1(\theta_1;\xi)$, and $\nu_2(\theta_2;\xi)$ such that $\nu_1(\theta_1;0)=u_L$, and $\nu_2(\theta_2;0)=u_R$ by minimizing the following loss functions ($k=1,2$) 
  \begin{eqnarray*}
    \mathcal{L}_k(\theta_k) & = & \|\nu'_k(\theta_k;\cdot) - r_k(\nu_k(\theta_k;\cdot))\|_{L^2(\R_+;\R^2)} \, ,
  \end{eqnarray*}
We denote by $\nu^*_k$ the optimized networks.
\item We numerically determine $\xi^{*} >0$ such that $\nu^*_1(\xi^{*})=\nu^*_2(\xi^{*})$ and corresponding to the intermediate state $u_1^*$.
\color{red}
 \end{enumerate} 
Once these intermediate waves identified, we can use them as initial condition for NDNN method.
\\
\noindent
{\bf Generation of 1 shock and 1 rarefaction wave.} 
Without loss of generality, let us assume that the solution to \eqref{ibvp2bis} is constituted by a 1-shock and a 2-rarefaction waves.  Keeping in mind that the domain decomposition which is proposed in this paper is only isolated regions between shock waves, we are only  interested in the evaluation of the value of $u_1^*$ such that $(u_L,u_1^*)$ is a 1-shock. In this goal we proceed as in the case of 2 shock waves, except that the searched intermediate state is now at the intersection of a 1-shock curve $s_1$, issued from $u_L$
\begin{eqnarray*}
  \left. 
  \begin{array}{rcl}
    s'_1(\xi) & = & r_1(s_1(\xi)), \, \textrm{ with } \xi > 0 \, ,\\
  s_1(0) & = & u_L \, ,
      \end{array}
  \right.
  \end{eqnarray*}
 and of the 2-rarefaction curve graph $w_2$,  and issued from $u_R$ 
 \begin{eqnarray*}
     \left. 
  \begin{array}{rcll}
  w'_2(\xi) & = & r_2(w_2(\xi)), &\textrm{ with } \xi < \lambda_2(u_R) \, ,\\
  w_2(\lambda_2(u_R)) & = & u_R \, . & 
     \end{array}
 \right.
 \end{eqnarray*}
 Finally we solve  $s_1(\xi^*)=w_2(\xi^*)$ and define $u_1^*=s_1(\xi^*)=w_2(\xi^*)$. The neural network-based  algorithm is similar as above.
%
% As above the algorithm reads as 
%\begin{enumerate}
%\item we optimize 2 vector-valued neural networks $\nu_1(\xi)$ (approximating the 1-shock curve), and $\nu_2(\xi)$ (approximating the 2-rarefaction curve) such that $\nu_1(0)=u_L$, and $\nu_2(\lambda_2(u_R))=u_R$ by minimizing the following loss functions ($k=1,2$) 
%  \begin{eqnarray*}
%\mathcal{L}_k(\bs{\Theta}) & = & \|\nu'_k(\cdot) - r_k(\nu_k(\cdot))\|_{L^2(\R_+;\R^2)}
%    \end{eqnarray*}
%\item We numerically determine $\overline{\xi}$ such that $\nu_1(\overline{\xi})=\nu_2(\overline{\xi})$.
 %  \end{enumerate}
 %
\subsubsection{Initial wave decomposition for 3 equation systems $(m=3,\, \mu=3$)}
{
 %In this section, we detail the wave decomposition for \eqref{pde+bc;d} for $m=3$.  
 The additional difficulty compared to the case $m=2$, is that the starting and ending states connected by the 2-wave are both unknown. We then propose the following approach.
\begin{enumerate}
\item We determine the 1st and 3rd wave respectively issued from $u_L$ and $u_R$. The corresponding curves are defined by $v_k$  and parameterized by real variables $\xi_k$ ($k=1,3$)
\begin{eqnarray*}
  v'_k(\xi_k) & = & r_k(v_k(\xi_k)) \, , 
  \end{eqnarray*}
with corresponding initial condition depending on the type of waves. For instance, $v_1(0)=u_L$ if the first curve is a 1-shock curve.
\item The intermediate curve $v_2$ parameterized by $\xi_2$ satisfies
\begin{eqnarray*}
  v'_2(\xi_2) & = & r_2(v_2(\xi_2)) \, , 
\end{eqnarray*}
and such that for some $\xi_{1}^{*}$, $\xi_{2}^{*}$, and $\xi_{3}^{*}$:
\begin{itemize}
\item $v_2(0) = v_1(\xi_{1}^{*})$ and $v_2(\xi_{2}^{*}) = v_3(\xi_{3}^{*})$ if the second curve is a 2-shock curve;
\item $v_2(\lambda_2(v_2(\xi_{2}^{*})) = v_3(\xi_{3}^{*})$ and $v_2(\xi_{2}^{*}) = v_3(\xi_{3}^{*})$ if the second curve is a 2-rarefaction curve.
\end{itemize}
The intermediate states then correspond to $u_1^{*}=v_2(\xi^{*}_2)$ and $u_2^{*}=v_3(\xi^{*}_3)$.
\end{enumerate}
Let us discuss the corresponding computational algorithm. Let $\nu_k$ denote the networks approximating $v_k$, an by $I_k$ the domain in $\xi$ (typically a real interval). 
We first minimize the following local loss functions for $k=1$ and $k=3$
\begin{eqnarray*}
\mathcal{L}_k(\theta) & = &\|\partial_{\xi_k} \nu_k(\theta_k;\cdot)-r_k(\nu_k(\theta_k;\cdot))\|_{L^2(I_k)}^2 
+
C_k\, 
\end{eqnarray*}
where $C_k$ denotes the corresponding initial condition, 
and denote by $\nu^*_1$ and $\nu^*_3$ the optimized networks.
Then say for a $2$-shock curve, $\nu_2^*$ is the network which optimizes the loss functional
\begin{eqnarray*}
\left.
\begin{array}{lcl}
\mathcal{L}_2(\theta_2) &= & \,\,\,\, \gamma_1\|\partial_{\xi_2} \nu_2(\theta_2;\cdot)-r_2(\nu_2(\theta_2;\cdot))\|_{L^2(I_2)}^2 + \gamma_2 \|\nu_2(\theta_2;0)-\nu^*_1(\xi_1)\|^2_{L^2(I_2)} \\
& & + \gamma_3 \|\nu_2(\theta_2;\xi_2)-\nu^*_3(\xi_3)\|^2_{L^2(I_2)} \, ,
\end{array}
\right.
\end{eqnarray*}
for some positive hyper-parameters $\gamma_1$, $\gamma_2$ and $\gamma_3$.  Once the wave decomposition is performed, the corresponding function can be taken as initial condition within the NDNN method.}

\color{black}
\section{Gradient descent algorithm and efficient implementation}\label{sec:ddm}
In this section we discuss the implementation of gradient descent algorithms for solving the
minimization problems  \eqref{e:min(L)}, \eqref{e:min(L),D} and \eqref{e:min(L),D;m}.
We note that these problems involve a global loss functional measuring the residue of
HCL in the whole domain, as well Rankine-Hugoniot conditions, 
which results in training of a number of neural networks.
In all the tests we have done, the gradient descent method converges and provides accurate results.
We note also, that in  problems with a large number of DLs, 
the global loss functional couples a large number of networks and the gradient descent algorithm may converge slowly. For these problems we present a domain decomposition method (DDM).

\subsection{Classical gradient descent algorithm for HCLs}\label{subsec:gda}
All the problems \eqref{e:min(L)}, \eqref{e:min(L),D} and \eqref{e:min(L),D;m} being similar, 
we will demonstrate in details the algorithm for the problem \eqref{e:min(L),D}.
We assume that the solution is initially constituted by 
i) $D\in\{1,2,\ldots,\}$ entropic shock waves emanating from $x_1,\ldots,x_D$, 
ii) an arbitrary number of rarefaction waves, and that 
iii) there is no shock generation for $t\in [0,T]$.

\begin{algorithm}
	\begin{algorithmic}[1]
		\STATE set $k=0$
		\REPEAT 
		\STATE  compute the gradient: $\nabla\bs{\cal L}(\bs{\Theta}^k)$
		\STATE  update the parameters: $\bs{\Theta}^{k+1} = \bs{\Theta}^k - \lambda_k \nabla\bs{\cal L}(\bs{\Theta}^k)$
		\STATE  set $k=k+1$ 
		\UNTIL{ \  $\sum_{i=1}^D\|{\bf n}_i^k-{\bf n}_i^{k-1}\|_{L^2(0,T)}<\epsilon$}
	\end{algorithmic}
	\caption{
		Gradient descent algorithm  \newline
		\textbf{Input:} \newline
			\hspace*{0.5cm} $\epsilon > 0$ --- tolerance \newline
		\hspace*{0.5cm} $\lambda_k > 0$ --- learning rate (sufficiently small) \newline
		\hspace*{0.5cm} $\bs{\Theta}^0\in\Theta_{ad}$ --- initial guess \newline
		\textbf{Output:} \newline
		\hspace*{0.5cm} $\bs{\Theta}^k$ --- approximation to the solution 
		$\bs{\Theta}^*$ of the problem 		\eqref{e:min(L),D}}
	\label{alg:gradient}
\end{algorithm}

The gradient algorithm applied to the problem \eqref{e:min(L),D} is given in Algorithm \ref{alg:gradient}.
We note that the functional $\bs{\cal L}(\bs{\uptheta})$ involves $L^2$ norms, which are not appropriate
in the meshless context of neural network computing.
Instead we approximate the integrals with sums over a number of points called  ``learning points".
More specifically, we choose 
a finite set of points ${\cal P}_R$ in the rectangle $R=[0,1]\times[0,T]$,
another finite set of points ${\cal P}_{(0,1)}$, and
another finite set of points ${\cal P}_{(0,T)}$.

Then the functional $\bs{\cal L}(\bs{\theta})$  in  \eqref{e:min(L),D} is approximated as follows
\begin{eqnarray}
\mathcal{L}(\bs{\Theta}) 
& = & 
\lambda 
\sum_{i=1}^{D+1}
\sum_{P\in{\cal P}_R}
\big|
(\partial_t {\bf N}_i(\cdot) 
+
\partial_x f({\bf N}_i(\cdot)))\circ {\bf T}_i(P)
\big|^2
\nonumber\\
&& +
\mu
\sum_{i=1}^D
\sum_{P\in{\cal P}_{(0,T)}}
\big|\partial_t {\bf n}_i(P)
\big[{\bf N}_i^+({\bf n}_i(P),P) - {\bf N}_i^-({\bf n}_i(P),P)\big] - 
\nonumber\\
&&
\hspace*{39mm}
\big[f({\bf N}_i^+({\bf n}_i(P),P))- f({\bf N}_i^-({\bf n}_i(P),P))
\big|^2
\nonumber\\
& & + 
\kappa
\sum_{i=1}^{D+1}
\sum_{P\in{\cal P}_{(0,1)}}
|{\bf N}_i^-({\bf T}_i(P,0)) - u_0({\bf T}_i(P,0))|^2.
\label{e:L,Dsw;sum}
\end{eqnarray}

The gradient descent algorithm associated to the minimization of \eqref{e:L,Dsw;sum}, which we use in the computations,
is given by Algorithm \ref{alg:gradient} with $\bs{\cal L}(\bs{\uptheta})$ given by \eqref{e:L,Dsw;sum}.
The global solution ${\bf N}^k$ at the iteration $k$ is then  constructed as  follows.
We denote by ${\bf N}_i^k$, resp. ${\bf n}_i^k$, the networks with parameters $\bs{\Theta}_i^k$, 
resp. $\bs{\uptheta}_i^k$,
and let $Q_i^k$ be the domain bounded by networks ${\bf n}_{i-1}^k$ and ${\bf n}_i^k$.
Then ${\bf N}^k(x,t) = {\bf N}_i^k({\bf T}_i^{-1}(x,t))$ with $(x,t)\in Q_i^k$.
%
%Notice that in the case of HSCL with $m$ equations, the minimization procedure  is actually almost identical as %\eqref{e:L,grad}, except that the $L^2$-norms are now performed in $L^2(0,T;\R^m)$ and $L^2(\Omega_0^i;\R^m)$.

\subsection{Gradient descent and domain decomposition methods}\label{subsec:ddm1}
Rather than minimizing the global loss function \eqref{e:L,Dsw} (or  \eqref{e:L,1sw},  \eqref{e:L,nsw2}), we here propose to decouple the optimization of the neural networks, and make it scalable. 
The approach is closely connected to domain decomposition methods (DDMs)
Schwarz Waveform Relaxation (SWR) methods \cite{gander2007optimized,HypSWR1,HypSWR2}. 
The resulting algorithm allows for embarrassingly parallel computation of minimization of local loss functions.

The algorithm is as follows.
For each $i=1,\ldots,D+1$, we introduce the networks
${\bf N}_i$ with parameters $\bs{\theta}_{i}$ and the two networks 
${\bf n}_i^l$, resp. ${\bf n}_i^r$, 
with ${\bf n}_i^l(0)=x_{i-1}$ and parameters $\bs{\uptheta}_i^l$,
resp.
with ${\bf n}_i^r(0)=x_i$ and parameters $\bs{\uptheta}_i^r$,
and set $\bs{\Theta}_i=\bs{\theta}_i\times \bs{\uptheta}_i^l\times \bs{\uptheta}_i^r$.
The networks ${\bf n}_i^l$ and ${\bf n}_i^r$ define a domain denoted by $Q_i$, 
see Section \ref{subsec:SW2} for notations. Like in Section \ref{subsec:SW2} we introduce by ${\bf T}_i:R\rightarrow Q_i$ the transformation defined by 
${\bf T}_i(x,t)=(({\bf n}_i^r(t)-{\bf n}_i^l(t))x+{\bf n}_i^l(t),t)$,
and for $i=1,\ldots,D+1$ consider the local functional 
\begin{eqnarray}
\bs{\mathcal{L}}_i(\bs{\Theta}_i) 
& = & 
\lambda 
\big\|
(\partial_t {\bf N}_i(\cdot) 
+
\partial_x f({\bf N}_i(\cdot)))\circ{\bf T}_i
\|^2_{L^2(R)} 
\nonumber\\
&+&
\mu
\Big(
\big\|\partial_t {\bf n}_i^l(\cdot)
\big[{\bf N}_i({\bf n}_i^l(\cdot),\cdot) - {\bf N}_{i-1}({\bf n}_i^l(\cdot),\cdot)\big] - \nonumber\\
&&
\hspace*{19mm}
\big[f({\bf N}_i({\bf n}_i^l(\cdot),\cdot))- f({\bf N}_{i-1}({\bf n}_i^l(\cdot),\cdot))\big]\|_{L^2(0,T)}^2
+
\nonumber\\
& & 
\hspace*{5mm}
\big\|\partial_t {\bf n}_i^r(\cdot)
\big[{\bf N}_{i+1}({\bf n}_i^r(\cdot),\cdot) - {\bf N}_i({\bf n}_i^r(\cdot),\cdot)\big] - \nonumber\\
&&
\hspace*{19mm}
\big[f({\bf N}_{i+1}({\bf n}_i^r(\cdot),\cdot))- f({\bf N}_i({\bf n}_i^r(\cdot),\cdot))\big]\|_{L^2(0,T)}^2
\Big)
\nonumber\\
&+ & 
\kappa
\|{\bf N}_i({\bf T}_i(\cdot,0)) - u_0({\bf T}_i(\cdot,0))\|^2_{L^2(0,1)},
\nonumber\\
&+&
\nu
\Big(
\|{\bf n}_i^l(\cdot)-{\bf n}_{i-1}^r(\cdot)\|_{L^2(0,T)}^2
+
\|{\bf n}_i^r(\cdot)-{\bf n}_{i+1}^l(\cdot)\|_{L^2(0,T)}^2
\Big),
\label{e:L,Dsw;DDM}
\end{eqnarray}
for some positive parameters $\lambda,\mu$, $\kappa$ and $\nu$, and where ${\bf N}_0=u(a)$, ${\bf N}_{D+2}=u(b)$.

We note that the minimization of $\bs{\cal L}_i$ in \eqref{e:L,Dsw;DDM} corresponds to the solution of the problem
\eqref{e:HCL,Qi;D} in the domain $Q_i$, coupled with a modified version of the Rankine-Hugoniot conditions \eqref{e:RH-1d;D}. This approach is very similar to a SWR methods with Rankine-Hugoniot transmission conditions and is referred hereafter as a ``Domain Decomposition method".

\begin{algorithm}
	\begin{algorithmic}[1]
		\STATE $\ell=0$
		\REPEAT 
			\FOR {$i=1,\ldots,D+1$}
				\STATE  minimize $\bs{\cal L}_i(\bs{\Theta}_i)$ with gradient descent algorithm; 
						let	$\bs{\Theta}_i^k$ be the solution
				\STATE  update $\bs{\Theta}_i=\bs{\Theta}_i^k$
			\ENDFOR
			\STATE  set $\ell=\ell+1$ 
			\STATE  set ${^\ell}\bs{\Theta}_i=\bs{\Theta}_i^k$ 
		\UNTIL{ \  $\sum_{i=1}^{D}\| {\bf n}_i^l - {\bf n}_{i-1}^r\|_{L^2(0,T)} <\delta^{\textrm{cvg}}$}
	\end{algorithmic}
	\caption{
		Domain decomposition algorithm  	\newline
		\textbf{Input:} \newline
		\hspace*{0.5cm} $\delta^{\textrm{cvg}} > 0$ --- tolerance \newline
		\hspace*{0.5cm} ${^0}\bs{\Theta}_i$ --- initial guesses of networks
		\newline
		\textbf{Output:} \newline
		\hspace*{0.5cm} ${^\ell}\bs{\Theta_i}$ --- approximation of the networks minimizing \eqref{e:L,Dsw;DDM}}
	\label{alg:ddm}
\end{algorithm}

DDM algorithm is summarized in Algorithm \ref{alg:ddm}. 
{Convergence is reached for some prescribed $\delta^{\textrm{cvg}}$. For instance in the framework of Schr\"odinger equation \cite{antoine2015analysis} solving with finite difference or element methods, $\delta^{\textrm{cvg}}=10^{-14}$. In this paper, we have however considered larger values  $10^{-5}$.}

We conclude this section by a discussion on the computational complexity of the NDNN vs DDM approaches.
Let us recall that in the NDNN method $(2D+1)$ neural networks are coupled within a global loss function $\bs{\mathcal{L}}$. 
We denote by $n_{u}$ (resp. $n_\gamma$) the number of parameters associated to ${\bf N}_i$ (resp. ${\bf n}_i$)
approximating the local solution in $Q_i$ (resp. the DLs $\gamma_i$). 
If $k_{\textrm{Global}}$ is the total number of the gradient method iterations to minimize 
$\bs{\mathcal{L}}$ up to a given tolerance, in $\R^{(D+1)n_u+Dn_{\gamma}}$, 
then the complexity of the direct method is $O(k_{\textrm{Global}} D(n_u+n_{\gamma}))$. Typically, this minimization is performed in parallel thanks to stochastic versions of gradient methods \cite{SGD}.

On the other hand, the DDM approach requires the optimization of the $(D+1$) neural networks ${\bf N}_i$, and of 
$D$ pairs of neural networks $({\bf n}_{i}^l,{\bf n}_i^r)$ approximating the DLs (notice that only $D$ may be sufficient but would require additional numerical study). 
We now denote by $k_{\textrm{Local}}$ the average number of gradient descent method iterations to minimize the {\it local} loss functions $\bs{\mathcal{L}}_i$ in $\R^{n_u+2n_{\gamma}}$ and by $k_{\textrm{DDM}}$ the average number of iterations to reach DDM convergence. The computational complexity is then given by $O(k_{\textrm{DDM}}k_{\textrm{Local}}D(n_u+n_{\gamma}))$. Recall that the DDM algorithm is trivially embarrassingly parallel. Moreover the local minimizations can also be performed with stochastic gradient methods providing hence a second level of parallelization. 

In conclusion, the DDM becomes relevant thanks to its scalability and for  $k_{\textrm{DDM}}k_{\textrm{Local}}<k_{\textrm{Global}}$, which  is expected for $D$ large.

\section{Numerics}\label{sec:numerics}

\subsection{Practical implementations}\label{subsec:numerics1}
This subsection is devoted to the practical aspects of the training process of neural networks. The implementation of the  algorithms above is performed using the library neural network ${\tt jax}$, see \cite{jax2018github}.
% \textcolor{red}
 {Although the algorithms look complex, they are actually very easy to implement using {\tt jax} and we did not face any difficulty in the tuning of the hyper-parameters. In this paper we propose a proof-of-concept of a novel method in low dimension, and which ultimately deals with simple (piecewise-)smooth functions. As a consequence, we have not addressed in details questions related to the choice of the optimization algorithm or of the hyper-parameters, because in this setting they are not particularly relevant. In our numerical simulations we have considered $\tanh$ neural networks with one or two hidden layers. The learning nodes to approximate the PDE residuals are randomly selected in the rectangular regions $R=(0,1)\times(0,T)$ (see Subsection \ref{subsec:SW1}). The weights $\lambda,\mu$ in \eqref{e:L,1sw} and \eqref{e:L,Dsw} are taken equal to $1/2$, and more generally for equations with several shock waves or for systems, an equal weight is given to each contribution of the loss functions.   Moreover the neural networks were designed to satisfy the initial condition and boundary conditions. In the tests below, the learning rate in the gradient descent method presented in Subsection \ref{subsec:gda}, is fixed to $2\times 10^{-3}$.  Some {\tt Python} codes are posted on {\tt github.com}.}

In all the numerical experiments below we consider the problem \eqref{e:HCL0;PDE}-\eqref{e:HCL;bc}, and in the following experiments we only specify $\Omega\times[0,T]$, $f(u)$ and $u_0$. We refer to the results with our algorithms as NDNN solution.

\subsection{Basic tests and convergence for 1 and 2 shock wave problems}\label{subsec:numerics2}
In this subsection, we do not consider any domain decomposition, so that only one global loss function is minimized as described in Subsections \ref{subsec:SW1}, \ref{subsec:SW2}.

\paragraph{\bf Experiment 1} 
In this  experiment we consider $\Omega\times[0,T]=(-4,1)\times[0,3/4]$ with $f(u)=4u(2-u)$. The initial data is given by 
\begin{eqnarray*}
  u_0(x)=
  \left\{
  \begin{array}{ll}
    1, & x < -2,\\
    \cfrac{1}{2}, & -2< x < 0,\\
    \cfrac{3}{2}, & 0<x.
    \end{array}
\right.
  \end{eqnarray*}
In the time interval $[0,1/2]$, it is constituted by a rarefaction and a shock wave with constant velocity. 
Then, in the time interval $[1/2,3/4]$ the initial shock wave interacts with the rarefaction wave to produce a new shock with non-constant velocity. More specifically the solution is given  by
\begin{eqnarray*}
 \underset{0\leq t\leq\frac{1}{2}}{u(x,t)}
 =
  \left\{
  \begin{array}{ll}
    1, & x < -2,\\
    1+ \cfrac{x+2}{8t}, & -2< x \leq -2+4t,\\
    \cfrac{1}{2}, & -2+4t<x < 0,\\
    \cfrac{3}{2}, & 0<x, 
    \end{array}
\right.
\hspace*{5mm}
\underset{\frac{1}{2}\leq t\leq\frac{3}{4}}{u(x,t)}
=
\left\{
\begin{array}{ll}
	1, & x < -2,\\
	1+ \cfrac{\gamma(t)-2}{8t}, & -2< x \leq \gamma(t), \\
	\cfrac{3}{2}, & \gamma(t)< x.
\end{array}
\right.
 \end{eqnarray*}
Here $\gamma$ is the DL and it solves
\begin{eqnarray*}
f\Big(\cfrac{3}{2}\Big) - f\Big(\cfrac{2+\gamma(t)}{8t}+1\Big) & = & \gamma'(t)\Big(\cfrac{1}{2} -\cfrac{\gamma(t)+2}{8t}\Big)\, ,
  \end{eqnarray*}
for $t\in [1/2,1]$ and $\gamma(1/2)=0$.

\begin{figure}[!ht]
  \begin{center}
    \hspace*{1mm}\includegraphics[height=3.8cm, keepaspectratio]{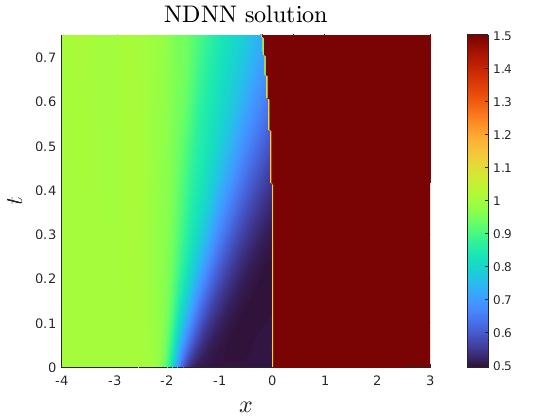}
      \hspace*{1mm}\includegraphics[height=3.8cm, keepaspectratio]{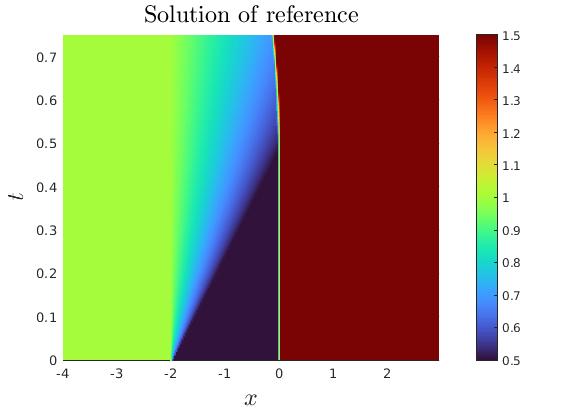}
     \hspace*{1mm}\includegraphics[height=3.8cm, keepaspectratio]{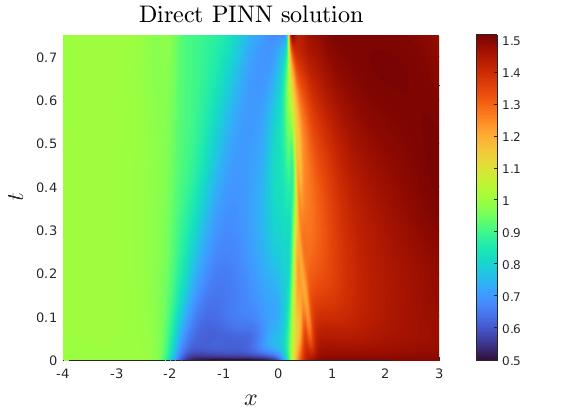} 
\caption{{\bf Experiment 1.} (Left) Neural network solution. (Middle) Solution of reference. (Right) Direct PINN solution.}
\label{figXa}
\end{center}
\end{figure}

\begin{figure}[!ht]
\begin{center}
    \hspace*{1mm}\includegraphics[height=8cm, keepaspectratio]{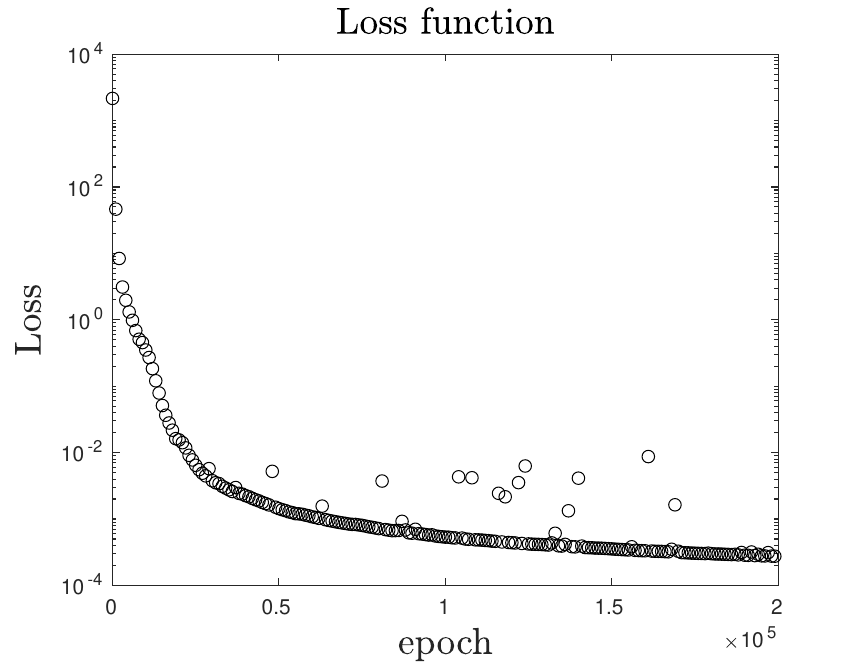}
\caption{{\bf Experiment 1.} Loss function.}
\label{figXb}
\end{center}
\end{figure}

 Numerically, we introduce 2 subdomains $\Omega^1_0=(-4,0)$ and $\Omega^2_0=(0,3)$. Notice that for the NN-algorithm, we use instead a regularized discontinuity for the non-entropic discontinuity located at $x=-1$. We introduce 3 neural networks -
 two space-time dependent and one time dependent, with 2 hidden layers and 20 neurons each and $2500$ learning nodes and the parameters. 
We report 
the neural network solution Fig. \ref{figXa} (Left),
the (exact) solution of reference Fig. \ref{figXa} (Middle),
the direct PINN solution Fig. \ref{figXa} (Right),
and the loss function as function of epoch number Fig. \ref{figXb}.

Let us mention that using the same numerical data, a direct PINN algorithm provides a very inaccurate approximation of the stationary then non-stationary shock waves,
while our algorithm provides accurate approximations. This last point is discussed in the 2 following tests.

\begin{figure}[!ht]
\begin{center}
  \hspace*{1mm}\includegraphics[height=6cm, keepaspectratio]{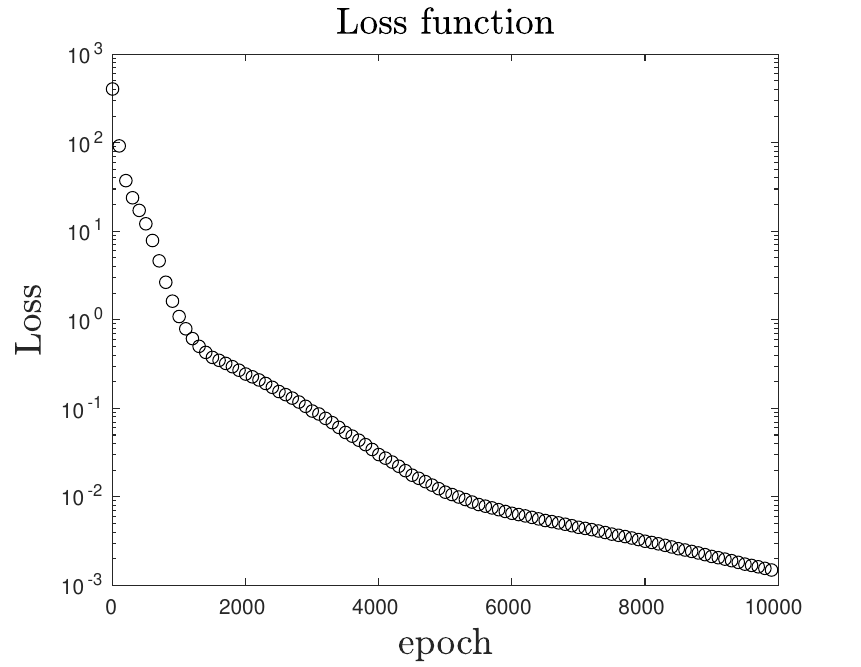}
    \hspace*{1mm}\includegraphics[height=6cm, keepaspectratio]{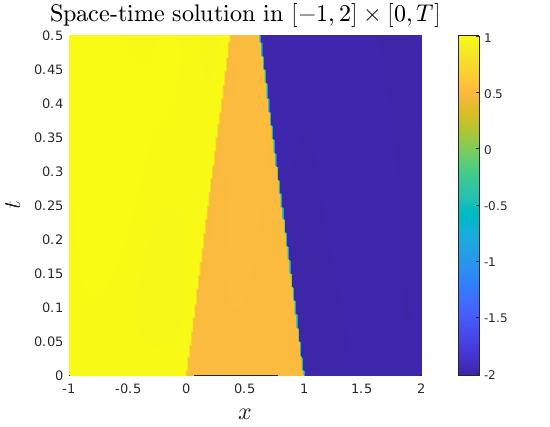}
\caption{{\bf Experiment 2.}(Left) Loss function. (Right) Space-time solution.}
\label{fig0}
\end{center}
\end{figure}

\begin{figure}[!ht]
\begin{center}
    \hspace*{1mm}\includegraphics[height=6cm, keepaspectratio]{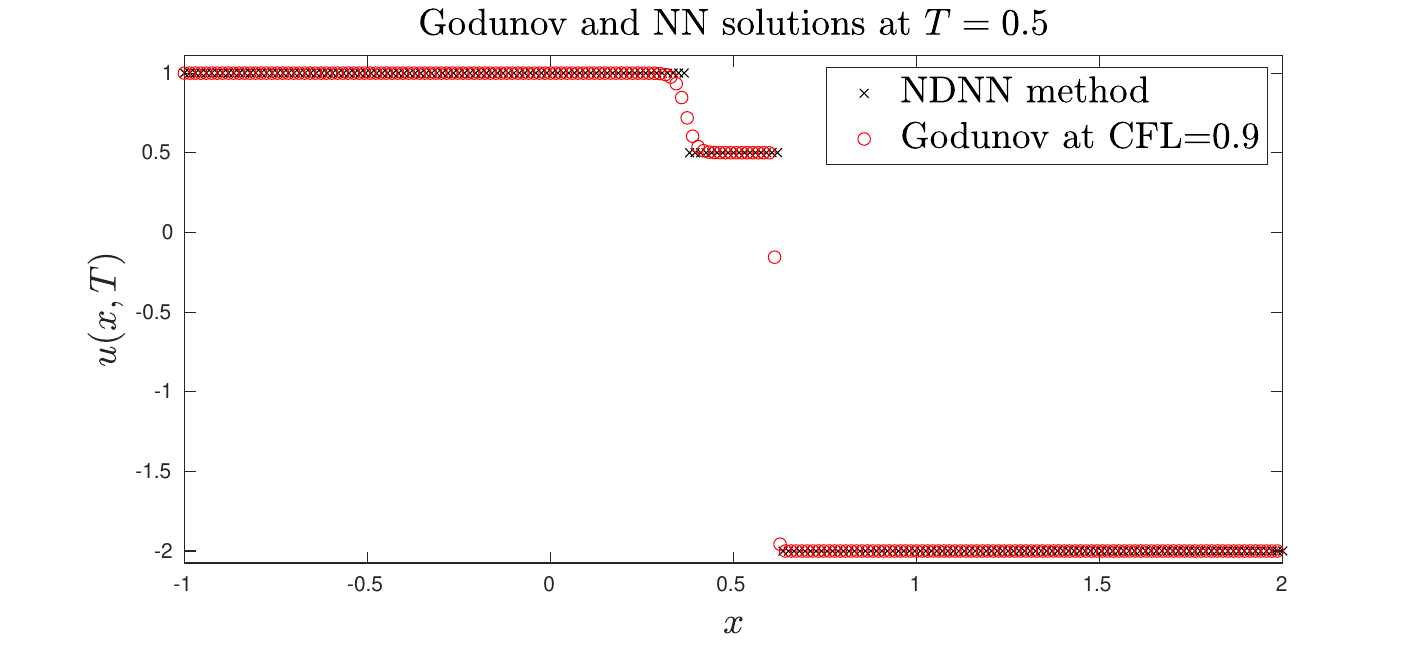}
\caption{{\bf Experiment 2.} (Left) Godunov scheme solution at CFL=0.9 and neural network solution at time $T=0.5$.}
\label{fig0b}
\end{center}
\end{figure}

\paragraph{\bf Experiment 2} 
Here $\Omega\times[0,T]=(-1,2)\times[0,0.5]$,  $f(u)=u^2/2$, and
\begin{eqnarray*}
	u_0(x)=
	\left\{
	\begin{array}{ll}
		1, & -1< x <0, \\
		\cfrac{1}{2}, & 0<x<1,\\
		-2, & 1<x <2.
	\end{array}
	\right.
\end{eqnarray*}
We compute the solution with the three following initial subdomains: $\Omega^1_0=(-1,0)$,  $\Omega^2_0=(0,1)$, $\Omega^3_0=(1,2)$ 
for $t\in [0,0.5]$. Note that the initial condition is constituted by two entropic shock waves. In this experiments the neural networks possess 1 hidden layer and 5 neurons each and $500$ learning nodes. In Fig. \ref{fig0} (Left) we report the loss function as function of epoch and in Fig. \ref{fig0} (Right), the space-time neural network solution at $T=0.5$ as well as the solution obtained by a Godunov scheme \cite{god1} with $200$ grid points at CFL=$0.9$ (CFL = $\Delta t \|f'(u)\|_{\infty}/\Delta x$).  We notice that unlike the Godunov scheme which naturally produces some numerical diffusion on both shock waves (particularly on the slowest one), the neural network solution is diffusion-free, see Fig. \ref{fig0b}. This is an interesting property which is more generally not shared with standard hyperbolic equation solvers.

\paragraph{\bf Experiment 3}
%\textcolor{red}
{In this experiment, we are specifically interested in the convergence of the NDNN algorithm. 
We consider a problem with 2 shock waves on $\Omega\times[0,T]=(-1,2)\times[0,2/5]$, with  $f(u)=u^2/2$ and
 \begin{eqnarray*}
  u_0(x)=
  \left\{
  \begin{array}{rr}
    2, & -1< x < 0,\\
    4x, & 0< x < 1,\\
    -4 , & 1<x<2.
    \end{array}
\right.
 \end{eqnarray*}
We solve the IVP for $t\in [0,2/5]$ in three different subdomains $\Omega^1_0=(-1,0)$, $\Omega^2_0=(0,1)$, $\Omega^3_0=(1,2)$. 
We use five networks ${\bf n}^1,{\bf n}^2,{\bf N}^1,{\bf N}^2,{\bf N}^3$ to approximate the DLs and the locals solutions, and each of the networks has one layer and the same number of neurons.
Each of the corresponding terms in the loss function has the same number of learning nodes. 
We approximate DLs amd the local solutions for different number of neurons and learning nodes.
 
We report the space-time solution on $(-1,2)\times [0,2/5]$ in Fig. \ref{figcv} (Left) as well as the total loss function in Fig. \ref{figcv} (Right). Regarding the convergence, we report $\ell^1-$norm errors of approximating the DLs over the time $[0,T]$ (denoted $\|\cdot\|_1$) as a function of neurons and of training$/$learning nodes. 
The error measures the $\ell_1$ norm of the different ${\bf n}^{i}-\gamma^i$, where 
${\bf n}^i$ is the approximation of DL and $\gamma^{i}$ is the DL of reference obtained with Godunov's method at CFL=0.99 and $1.5\times 10^{4}$ grid points. 
The error as a function of the number of neurons ($2$, $4$, $8$, $16$) and $512$ learning nodes is given in Fig. \ref{figcvB} (Left), 
and the error as a function of the the number of learning nodes ($8$, $18$, $72$, $288$, $1800$) 
with $16$ neurons is given in Fig. \ref{figcvB} (Middle). 
For the sake of completeness, we also report in Fig. \ref{figcvB} (Right) the graph of ${\bf n}^{i}$ 
with one-layer network and $8$ and $512$ learning nodes and $\gamma^{i}$ (lines of reference), the latter computed with Godunov method.
\begin{figure}[!ht]
\begin{center}
  \hspace*{1mm}\includegraphics[height=5cm, keepaspectratio]{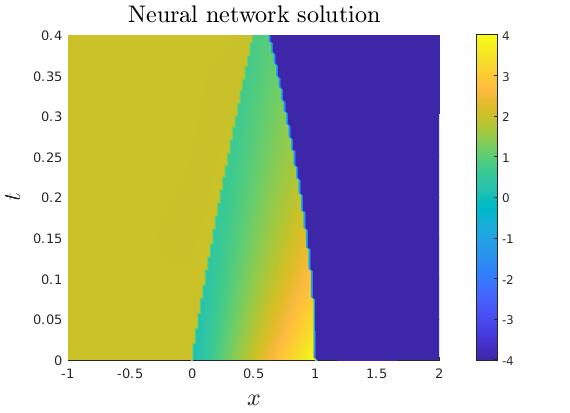}
   \hspace*{1mm}\includegraphics[height=5cm, keepaspectratio]{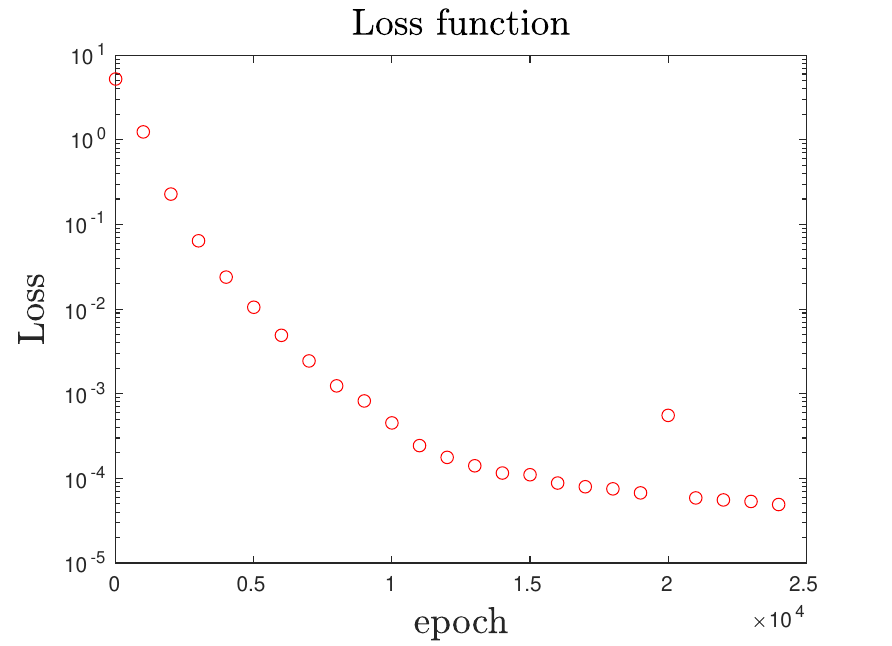}
\caption{{\bf Experiment 3.} (Left) Space-time solution (number of neurons is equal to $5$ with learning nodes is $1800$). (Right) Loss function.}
\label{figcv}
\end{center}
\end{figure}

\begin{figure}[!ht]
\begin{center}
  \hspace*{1mm}\includegraphics[height=3.8cm, keepaspectratio]{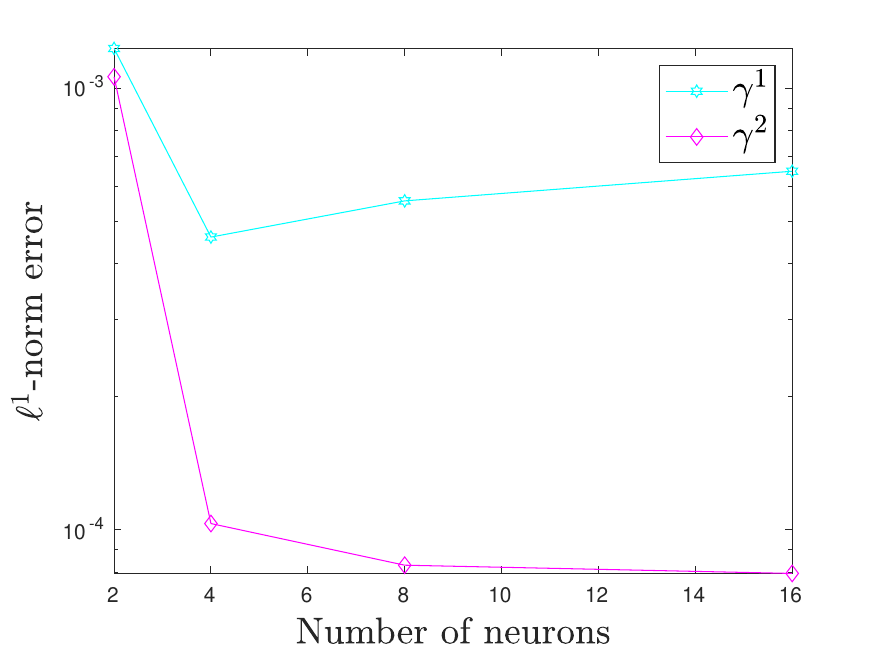}
  \hspace*{1mm}\includegraphics[height=3.8cm, keepaspectratio]{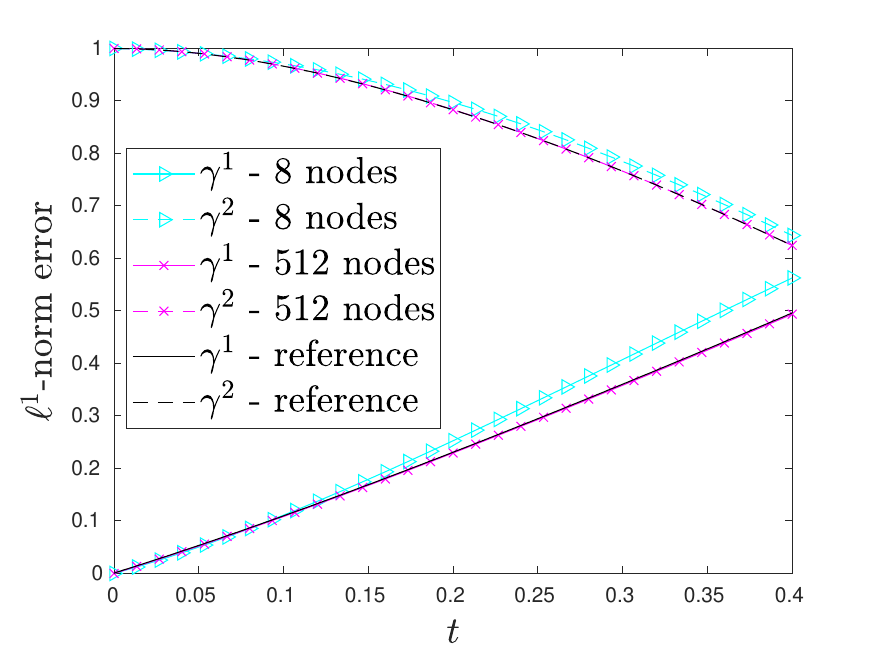}
   \hspace*{1mm}\includegraphics[height=3.8cm, keepaspectratio]{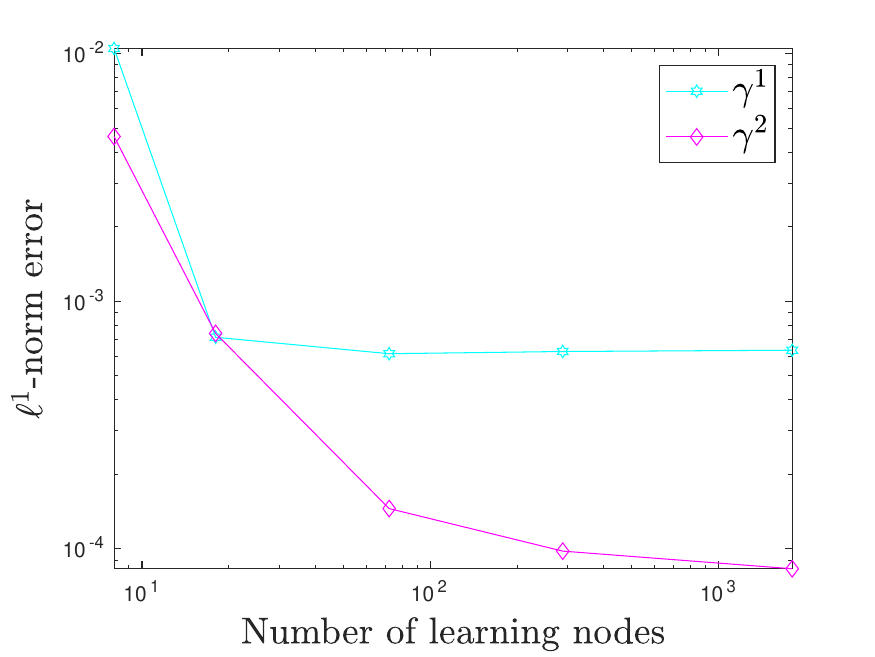}
\caption{{\bf Experiment 3.} (Left) $\ell^1$-norm in time: $\|{\bf n}^{i} -\gamma^{i}\|_{1}$ 
%\sout{for one-layer networks} 
with respectively $2$, $4$, $8$, $16$ neurons. (Middle) $\ell^1$-norm in time: $\|{\bf n}^{i} -\gamma^{i}\|_{1}$ 
%\sout{for one-layer networks} 
with $16$ neurons and respectively  $8$, $18$, $72$, $288$, $1800$ learning nodes. }
\label{figcvB}
\end{center}
\end{figure}
%\textcolor{red}
{These experiments allow to validate the convergence of the proposed approach.}
%\textcolor{blue}{tous les reseaux local-solutions et DL ont 1 layer avec le meme nombre de neurones et meme nombre de learning nodes. Le rafinement en learning nodes et neurones est donc applique a *tous* les reseaux, pas seulement sur les DL. Donc tout est bien consistent en terme de convergence.}

\subsection{Shock wave generation}\label{subsec:numerics3}
In this section, we demonstrate the potential of our algorithms to handle shock wave generation, as described in Subsection \ref{subsec:SW3}.
One of the strengths of the proposed algorithm is that it does not require to know the initial position$\&$time of birth, in order to accurately track the DLs. 
Recall that the principle is to assume that in a given (sub)domain and from a smooth function a shock wave will eventually be generated. 
Hence we decompose the corresponding (sub)domain in two subdomains and consider three neural networks: 
two neural networks will approximate the solution in each subdomain, and one neural network will approximate the DL. As long as the shock wave is not generated (say for $t<t^*$), the global solution remains smooth and the Rankine-Hugoniot condition is trivially satisfied (null jump); hence the DL for $t<t^*$ does not have any meaning.
\paragraph{\bf Experiment 4} 
We again consider the inviscid Burgers' equation, $\Omega\times[0,T]=(-1,2)\times[0,0.5]$  and
the initial condition
\begin{eqnarray*}
	u_0(x)=
	\left\{
	\begin{array}{ll}
		\cfrac{3}{4}-\tanh(2x), & -1< x < \cfrac{1}{2}, \\
		\cfrac{3}{4}-\tanh(2x), & \cfrac{1}{2}< x < \cfrac{3}{2},\\
		-\cfrac{1}{2}, & \cfrac{3}{2}< x< 2.
	\end{array}
	\right.
\end{eqnarray*}
Three initial subdomains are $\Omega^1_0=(-1,1/2)$,  $\Omega^2_0=(1/2,3/2)$, $\Omega^3_0=(3/2,2)$. Initially, an entropic  shock wave is located at $x=3/2$ at $t=0$.
The solution is smooth in the region covered by characteristics emanating from $\Omega^1_0\cup\Omega^2_0$ for $t<t^*$. 
Then a shock wave is generated at $t=t^*=3/5$. Hence, for $t<t^*$ the solution is constituted by one shock wave, and for $t>t^*$ by two shock waves. 

The solution in each subdomain is approximated by space-time neural networks with $30$ neurons and one hidden layer each and $900$ learning nodes.
We consider two time-dependent neural networks, one for approximating the chock wave initiated at $x=3/2$ and  the other initiated at $x=1/2$ to approximate the shock wave expected to be generated at time $t^*$.
The latter network has no particular meaning for $t<t^*$ - it separates the subdomains $\Omega^1_0$ and $\Omega^2_0$,
and models an artificial DL.

In Fig. \ref{fig1} (Left) we report the loss function as function of the epoch number.
In Fig. \ref{fig1} (Right) we report the corresponding reconstructed neural network space-time solution in the three subdomains by following the algorithm as explained in Subsection \ref{subsec:SW3}.
We observe the shock wave initiated at $t=0$,
and the generation of a shock wave between the subdomains $\Omega^1_t$ and $\Omega^2_t$ at $t^*$. 

We also report the neural network solution at $T=0$ and $T=3/5$ in Fig. \ref{fig2} (Left) as well as the graph of the neural network approximating the first and second DLs in Fig. \ref{fig2} (Middle). 
Finally, we report in Fig. \ref{fig2} (Right) the graph of the approximate flux jumps along the 
DLs as a function of time: $t \mapsto f(u(\gamma_{i}(t)^+,t))-f(u(\gamma_{i}(t)^-,t))$. We observe that along the first DL, the jump is close to zero until $t=t^*\approx 0.1$. This illustrates the fact that before $t<t^*$, there is no actual discontinuity along $\gamma_1$ (artificial discontinuity). For larger $t$, a jump appears in the flux (then on the solution) corresponding the generation of a shock wave. The second jump has a constant value as a function of $t$, which is consistent with the existence of a shock wave with constant velocity. 

\begin{figure}[!ht]
\begin{center}
  \hspace*{1mm}\includegraphics[height=6cm, keepaspectratio]{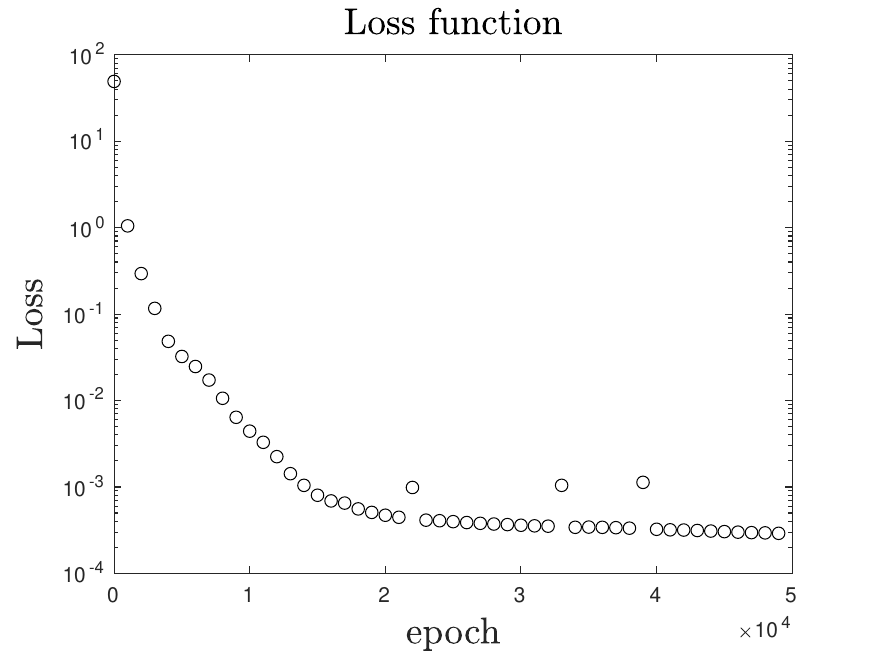}
    \hspace*{1mm}\includegraphics[height=6cm, keepaspectratio]{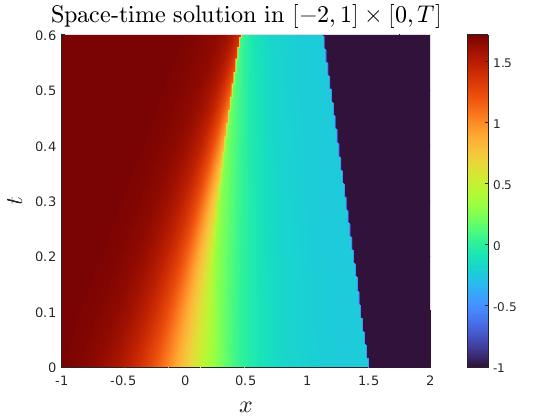}
\caption{{\bf Experiment 4.} (Left) Loss function. (Right) Space-time solution}
\label{fig1}
\end{center}
\end{figure}

\begin{figure}[!ht]
\begin{center}
  \hspace*{1mm}\includegraphics[height=3.8cm, keepaspectratio]{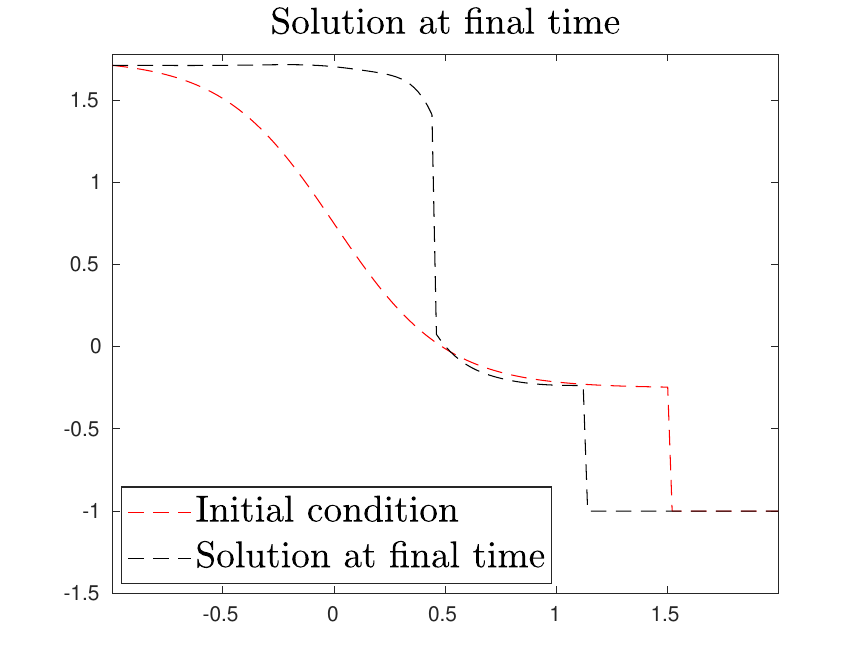}
  \hspace*{1mm}\includegraphics[height=3.8cm, keepaspectratio]{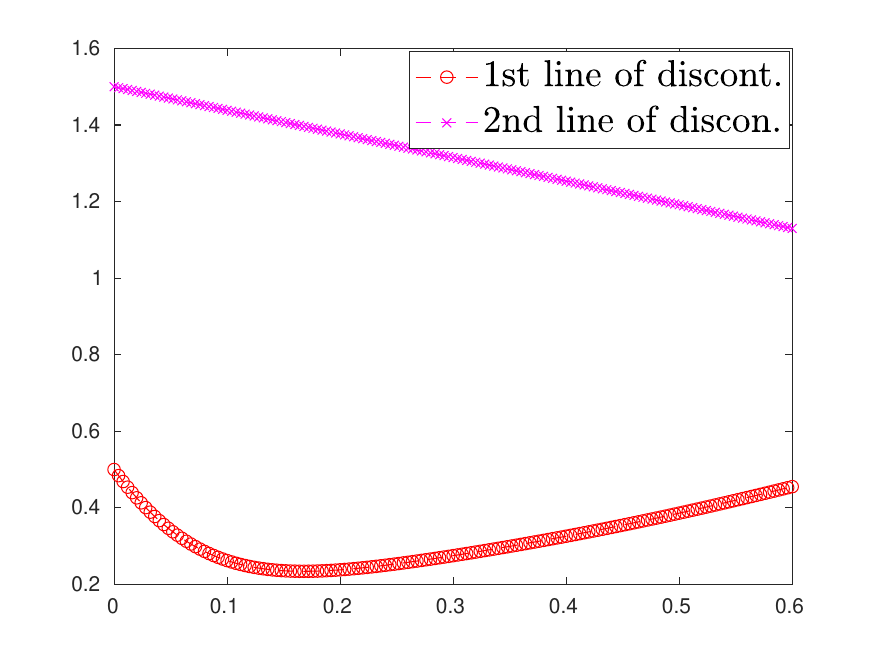}
    \hspace*{1mm}\includegraphics[height=3.8cm, keepaspectratio]{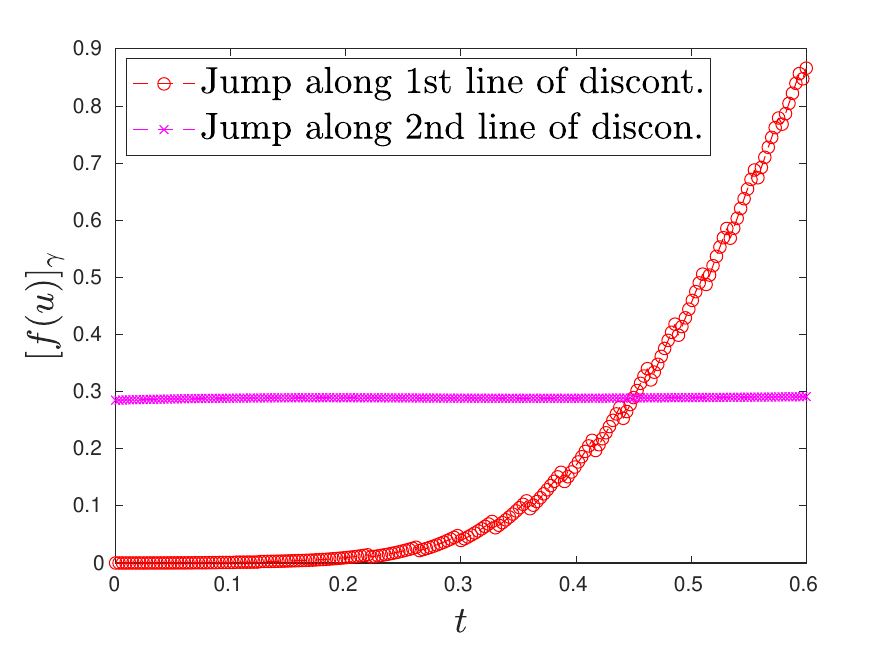}
\caption{{\bf Experiment 4.} (Left) Graph of the solution at $T=3/5$.  (Middle) Discontinuity lines. (Right) Flux jump along the DLs.}
\label{fig2}
\end{center}
\end{figure}

\subsection{Shock-Shock interaction}\label{subsec:numerics4}
In this subsection, we are proposing a test involving the interaction of two shock waves merging to generate a third shock wave. As explained in Subsection \ref{subsec:SW4}, in this case it is necessary re-decompose the full domain once the two shock waves have interacted. 
\paragraph{\bf Experiment 5}
% $f(u)=u^2/2$ and
% \begin{eqnarray*}
%  u_0(x)=
%  \left\{
%  \begin{array}{rr}
%    2, & -1< x < 0,\\
%    4x, & 0< x < 1,\\
%%%    -4 , & 1<x<2.
%    \end{array}
%\right.
%  \end{eqnarray*}
 The setting is identical to Experiment 3., with now $T=3/5$. We first solve the IVP for $t\in [0,3/5]$ in three different subdomains $\Omega^1_0=(-1,0)$, $\Omega^2_0=(0,1)$, $\Omega^3_0=(1,2)$. 
We report the space-time solution on $(-1,2)\times [0,3/5]$ in Fig. \ref{figI} (Left).  The time $t^*$ of interaction of the shocks is such that $\gamma_1(t^*)=\gamma_2(t^*)$. Numerically the ``exact'' time of shock wave interaction is computed by solving 
${\bf n}^1(t^*)={\bf n}^2(t^*)$, and which is numerically {\it a posterio} estimated at $t^*=0.45$ and located at $x^*=0.55$.   As a consequence $\Omega^2_{t^*}$ is reduced to one point. Beyond $t>0.45$, there is only one shock left, and the full space domain is now re-decomposed in two subdomains $\Omega^1_{t^*}=(-1,x^*)$ and $\Omega^2_{t^*}=(0.55,2)$ with a new DL located at $(x^*,t^*) \approx (0.55, 0.45)$. 
For $t>t^*$, we apply the same approach as before with two subdomains only. The general space-time solution can hence be reconstructed as presented in Fig. \ref{figI} (Right).
\begin{figure}[!ht]
\begin{center}
  \hspace*{1mm}\includegraphics[height=5cm, keepaspectratio]{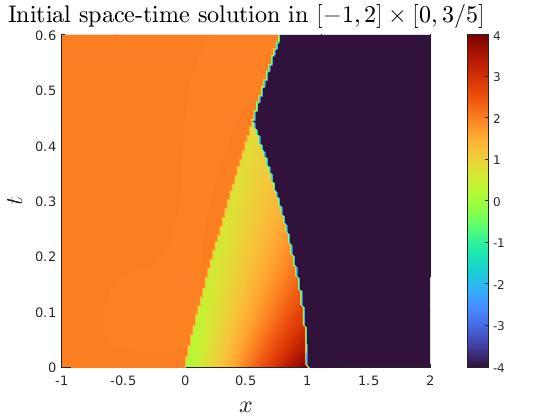}
   \hspace*{1mm}\includegraphics[height=5cm, keepaspectratio]{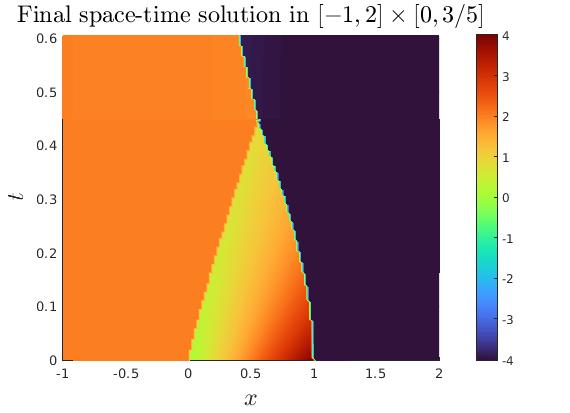}
\caption{{\bf Experiment 5.} 
	(Left) Space-time solution without shock interaction (artificial for $t>t^*=0.45$. 
	(Right) Space-time solution with shock interaction.}
\label{figI}
\end{center}
\end{figure}

\subsection{Entropy solution}\label{subsec:numerics5}
%\textcolor{red}
{We propose here an experiment dedicated to the computation of the viscous shock profiles and rarefaction waves and illustrating the discussion from Subsection \ref{subsec:SW5}. In this example, a regularized non-entropic shock is shown to be ``destabilized'' into rarefaction wave by the direct PINN method.}
\paragraph{\bf Experiment 6} We consider \eqref{e:HCL0} with $f(u)=u^2/2$, $\Omega_0=(-3,3)$ and
\begin{eqnarray*}
u_0(x) & = & \pm 2\Big(1 + \cfrac{-e^{4x}-0.01e^{-4x}}{e^{4x}+100e^{-4x}}\Big).
\end{eqnarray*}
We denote $u_{-} = u_0(-3^+) \approx 0$ and  $u_{+} = u_0(3^-) \approx -2$ 
(resp. $u_{-} \approx0$ and $u_{+} \approx 2$).  
The solution is obtained with networks with one hidden layer and $40$ neurons and $1600$ learning nodes per neural network. 
In Fig. \ref{fig3} we report the solution at time $T=0.65$, with two distinct initial data. 
As expected a viscous shock profile is captured when $f'(u_{-})>f'(u_{+})$ (resp. a rarefaction when $f'(u_{+})>f'(u_{-})$),  
which illustrates the entropic-like feature of {\it direct} PINN solvers, which will naturally be satisfied by our own neural network algorithm.
\begin{figure}[!ht]
\begin{center}
  \hspace*{1mm}\includegraphics[height=5cm, keepaspectratio]{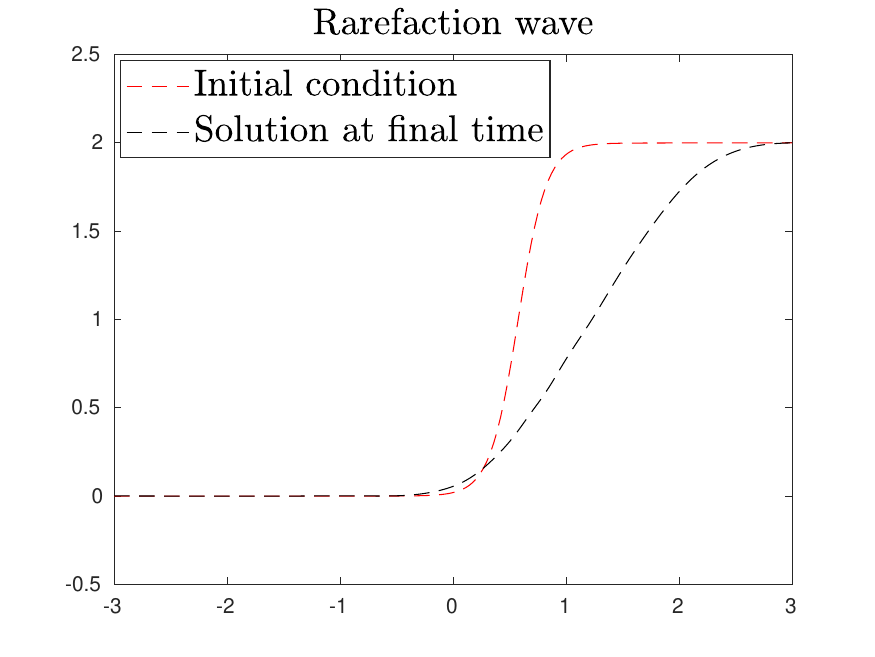}
    \hspace*{1mm}\includegraphics[height=5cm, keepaspectratio]{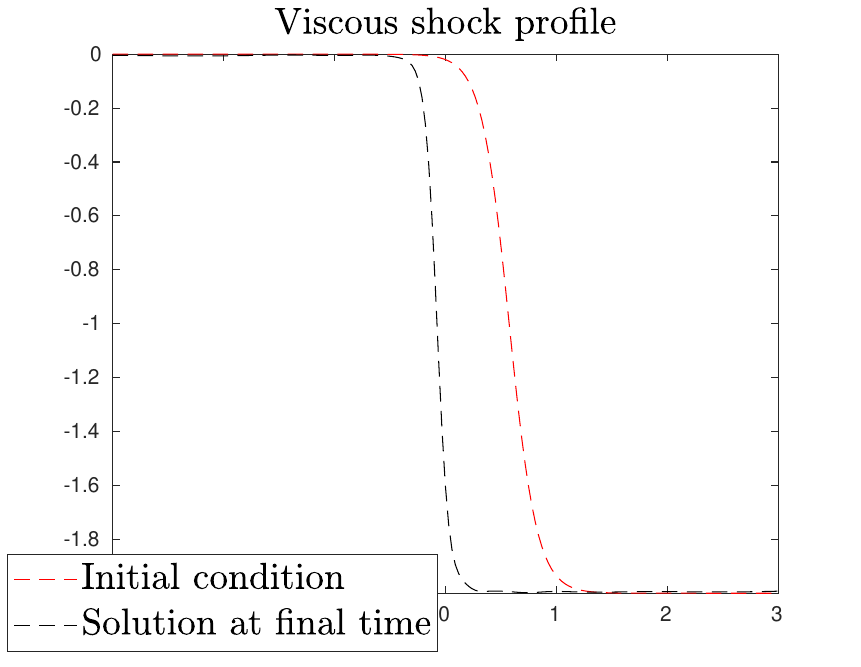}
\caption{{\bf Experiment 6.} (Left) Rarefaction for $f'(u_{+})>f'(u_{-})$. (Right) Viscous shock profile for $f'(u_{-})>f'(u_{+})$.}
\label{fig3}
\end{center}
\end{figure}

\subsection{Domain decomposition}\label{subsec:numerics6}
In this subsection, we propose an experiment illustrating the combination of the neural network based HCL solver developed in this paper with the domain decomposition method from Subsection \ref{subsec:ddm1}. We numerically illustrate the convergence of the algorithm. The neural networks have $30$ neurons and one hidden layer, and the number of learning nodes is $900$ learning nodes.
\paragraph{\bf Experiment 7} 
We consider $f(u)=u^2/2$ on $\Omega_0\times[0,T]=(-1,2)\times [0,1]$ with the following initial data
\begin{eqnarray*}
  u_0(x)=
  \left\{
  \begin{array}{rr}
    1, & -1<x < 0,\\
    2x & 0< x<1, \\
    0, & 1<x <2.
    \end{array}
\right.
  \end{eqnarray*}
We decompose $\Omega_0$ in three subdomains $\Omega^1_0=(-1,0)$,  $\Omega^2_0=(0,1)$,  $\Omega^3_0=(1,2)$. We implement the algorithm derived in Subsection \ref{subsec:ddm1}. 
We report the reconstructed solution at (Schwarz) convergence (after $50$ Schwarz iterations) in Fig. \ref{fig5}, and the solution at final time $T=1$ in Fig. \ref{fig6} (Left) and the local loss function values after $\ell_{\infty}=10^{5}$ optimization iterations for each Schwarz iteration  in Fig. \ref{fig6} (Right): 
that is $\bs{\cal L}_{i}(\bs{\Theta}_i^\ell)$, $i=1,2,3$, after $\ell$ optimization iterations. We observe that the local loss function values after a fixed number of optimization iterations ($\ell$) are decreasing as a function of 
(Schwarz) DDM iteration, illustrating the overall convergence of the Schwarz DDM algorithm.
\begin{figure}[!ht]
\begin{center}
    \hspace*{1mm}\includegraphics[height=8cm, keepaspectratio]{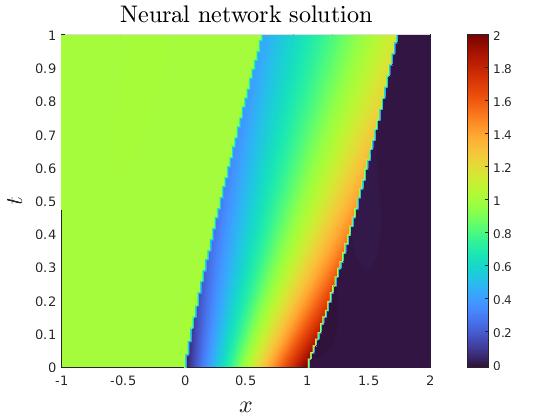}
\caption{{\bf Experiment 7.} Reconstructed space-time solution.}
\label{fig5}
\end{center}
\end{figure}

\begin{figure}[!ht]
\begin{center}
  \hspace*{1mm}\includegraphics[height=5cm, keepaspectratio]{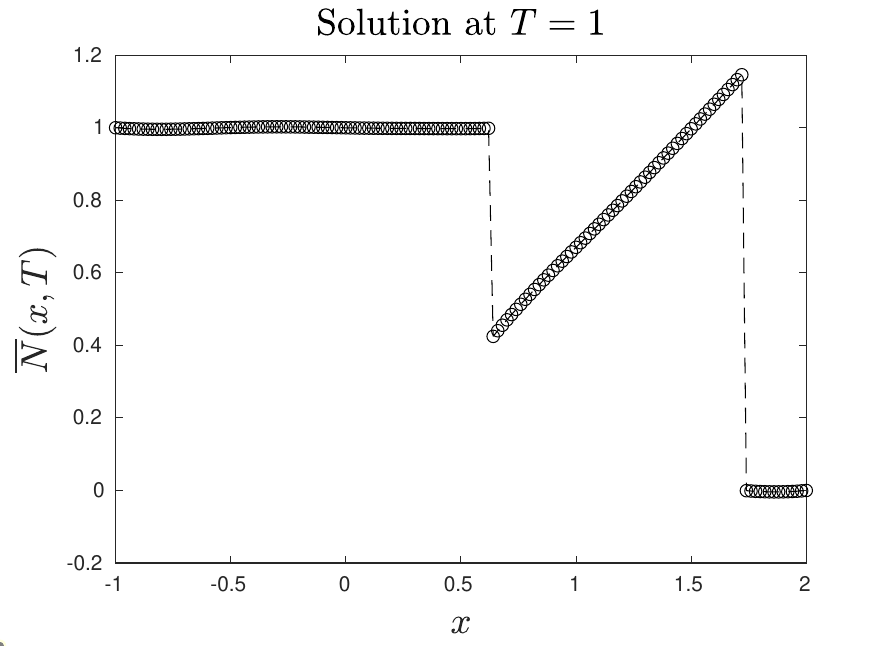}
    \hspace*{1mm}\includegraphics[height=5cm, keepaspectratio]{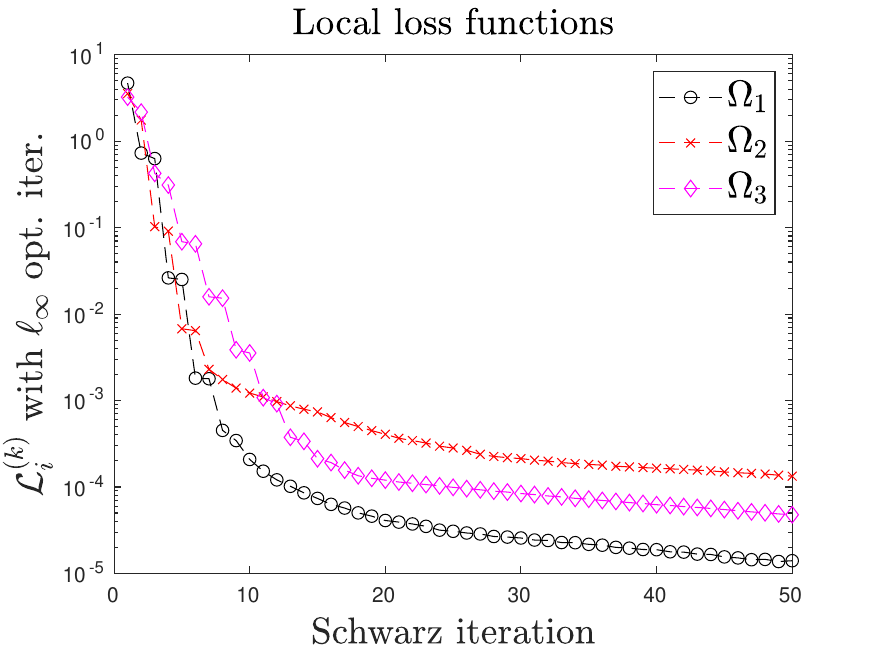}
\caption{{\bf Experiment 7.} (Left) Solution at $T=1$. (Right) Local loss function values after a fixed  number $\ell_{\infty}$ of optimization iterations.}
\label{fig6}
\end{center}
\end{figure}
The DDM naturally makes sense for much more computationally complex problems. This test however illustrates a proof-of-concept of the SWR approach.
%\input hypsysv18
% start
\subsection{Nonlinear systems}\label{subsec:numerics7}
In this subsection, we are interested in the numerical approximation of hyperbolic systems with shock waves.
\paragraph{\bf Experiment 8} In this experiment we focus on the initial wave decomposition for a Riemann problem.  The system considered here is the Shallow water equations ($m=2$). 
\begin{eqnarray}\label{sw1}
	\left.
	\begin{array}{lcl}
		\partial_t h + \partial_x(h u) & = & 0,\\
		\partial_t (h u) + \partial_x\Big(h u^2 +\cfrac{1}{2}gh^2\Big) & = & 0 \, ,
	\end{array}
	\right.
\end{eqnarray}
where $h$ is the height of a compressible fluid, $u$ its velocity, and  $g$ is the gravitational constant taken here equal to $1$.  The spatial domain is $(-0.1,0.1)$, the final time is $T=0.0025$, and we impose null Dirichlet boundary conditions. 
\\
%\textcolor{red}
{\noindent{\bf Experiment 8a.} 
The initial data is given by
\begin{eqnarray}\label{sw2}
	(h_0,h_0u_0)=
	\left\{ 
	\begin{array}{ll}
		(3,5), & x<0 \, , \\
		(3,-5), & 0<x \, .\\
	\end{array}
	\right.
\end{eqnarray}
Note that the corresponding solution is constituted by 2 entropic shock waves. 
We first implement the initial wave decomposition using the method proposed in Subsection \ref{subsec:IWD} with neural networks with 2 hidden layers and 5 neurons each, and $150$ learning nodes. The domain in $\xi$ is $[0,3/2]$. We report in Fig. \ref{figHL} (Left) the 1-shock and 2-shock curves (Hugoniot loci), that is  $\{(\xi,\nu_k(\xi)) \, : \, \xi \geq 0\}$. The loss functions are reported in \eqref{figHL} (Right).}
\begin{figure}[hbt!]
\begin{center}
  \includegraphics[height=4cm,keepaspectratio]{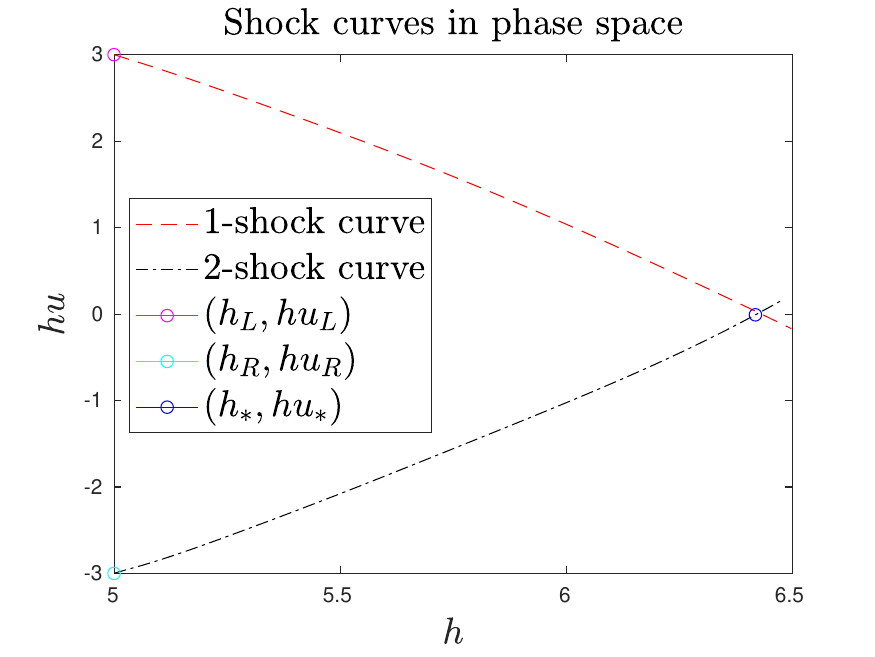}
  \includegraphics[height=4cm,keepaspectratio]{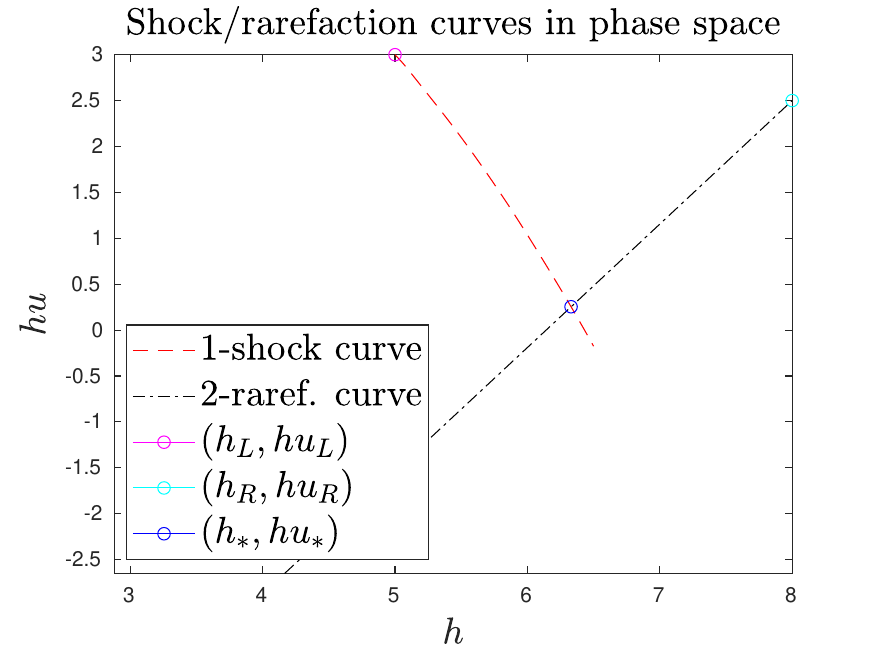}
    \includegraphics[height=4cm,keepaspectratio]{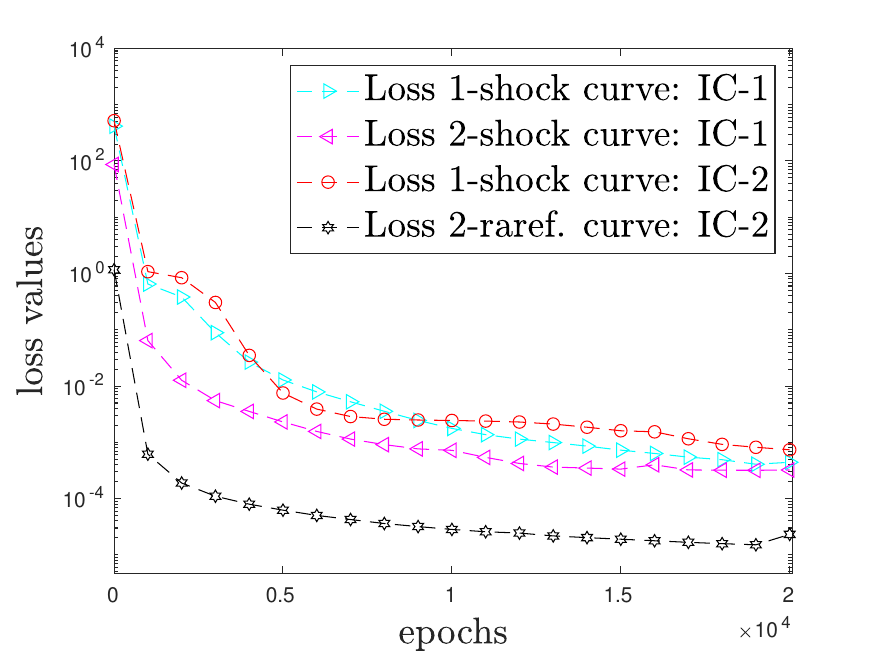}
\end{center}
\caption{{\bf Experiment 8a.} (Left) Shock curves in phase space with initial condition \eqref{sw2}. (Middle) 1-shock curve and 2-rarefaction curve with initial condition \eqref{sw3}. (Right) Loss functions for constructing simple waves \eqref{sw2} and \eqref{sw3}}
\label{figHL}
\end{figure}
%\textcolor{red}
{We numerically obtain $h_1^{*} \approx 6.428$ and $(hu)_1^{*} \approx 0.009$ and validate Lax entropy conditions for $1$-shock and $2$-shock waves $\lambda_k(u_R)<\sigma_k=\lambda_k(u_k^*)<\lambda_k(u_L)$ for $k=1,2$. This wave decomposition allows to define a new IVP at $t=0^{+}$, where the discontinuity between $(h_L,hu_L)$ and $(h_1^*,(hu)^*_1)$ represents a 1-shock, and the one between  $(h_1^*,(hu)^*_1)$ and $(h_L,hu_L)$ represents a 2-shock.}

%\textcolor{red}
{
For the sake of completeness we propose another decomposition generating a 1-shock and a 2-rarefaction wave. That is we consider
\begin{eqnarray}\label{sw3}
  (h_0,h_0u_0)=
  \left\{ 
  \begin{array}{ll}
(5,3), & x<0 \, , \\
(8,5/2), & 0<x \, .\\
    \end{array}
  \right.
\end{eqnarray}
 We apply the same method as above with 2 neural networks approximating $s_1$ and $w_2$ which are reported in Fig. \ref{figHL} (Middle). The intermediate state $u_1^*=(h^*_1,(hu)^*_1)$  involved in the 1-shock is finally numerically estimated as $h^*_1 \approx 6.3298$ and $(hu)^*_1 \approx 0.2543$. We also report the loss functions \eqref{figHL} (Right). Once the decomposition is done, it is then possible to apply the NDNN method.}
\\
\noindent{\bf Experiment 8b.}
Here we consider the NDNN method applied to the Shallow water problem \eqref{sw1} with initial condition \eqref{sw2}. We first perform the initial wave decomposition performed from Experiment 8a. then apply our NDNN method developed in Section \ref{sec:nd-solver}. We consider $8$ neural networks (approximating $h$, $hu$ and the two lines of discontinuity) each with $2$ hidden layers, $4$ neurons per layer. 
We report the space-time solution $h$, $hu$ in Fig. \ref{hyp_sys_init0a}.
\begin{figure}[!ht]
	\begin{center}
		\hspace*{1mm}\includegraphics[height=5.5cm, keepaspectratio]{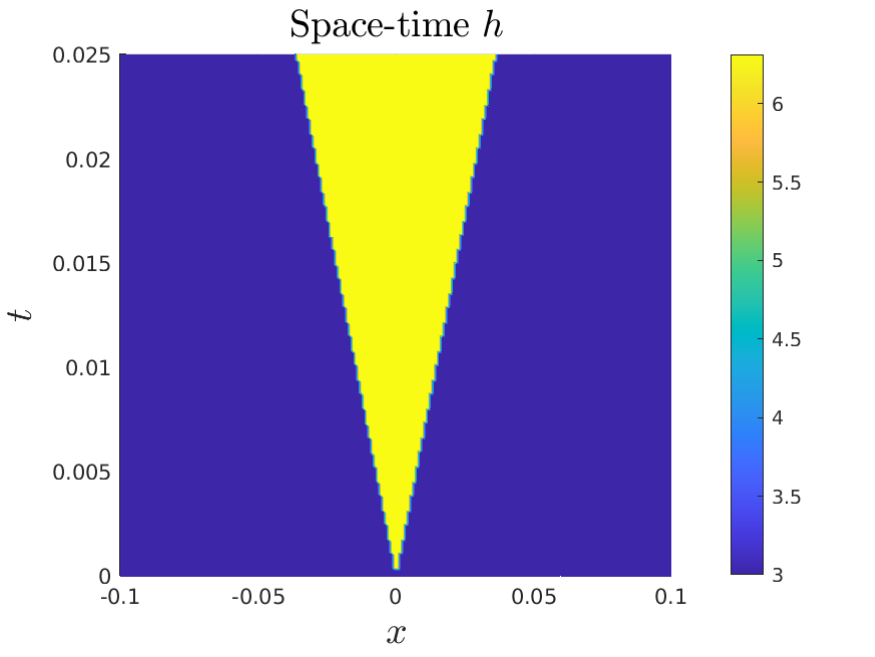}
		\hspace*{1mm}\includegraphics[height=5.5cm, keepaspectratio]{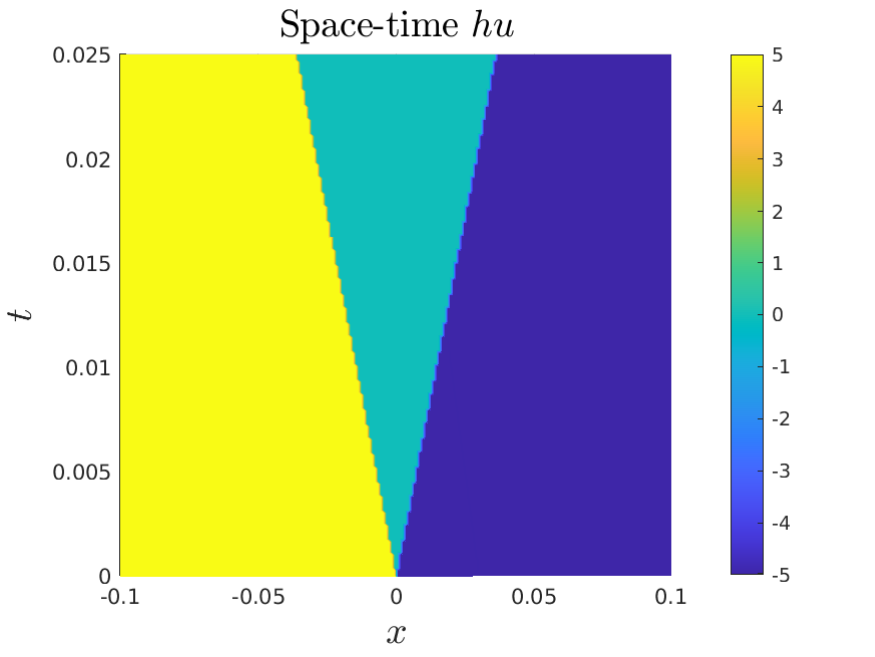}
		\caption{{\bf Experiment 8b.} Approximate space-time solution (Left)  $h:(x,t) \mapsto h(x,t)$. (Right) $hu: (x,t) \mapsto hu(x,t)$.}
		\label{hyp_sys_init0a}
	\end{center}
\end{figure}
We also report in Fig. \ref{hyp_sys_init0b} the initial and final approximate and exact solutions $h$ (Left) and $hu$ (Middle), as well as the loss function Fig. \ref{hyp_sys_init0b} (Right).
\begin{figure}[!ht]
	\begin{center}
		\hspace*{1mm}\includegraphics[height=3.8cm, keepaspectratio]{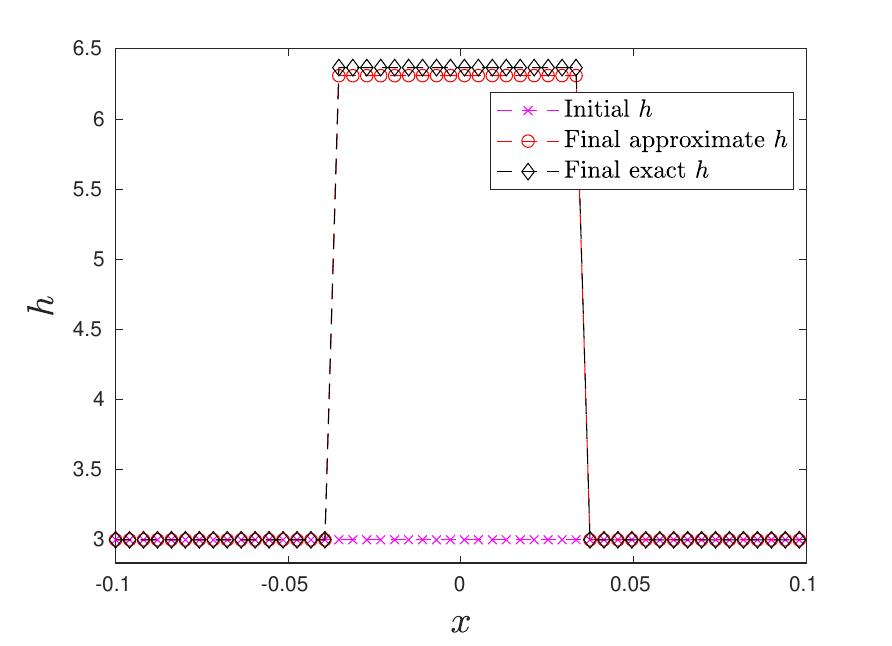}
		\hspace*{1mm}\includegraphics[height=3.8cm, keepaspectratio]{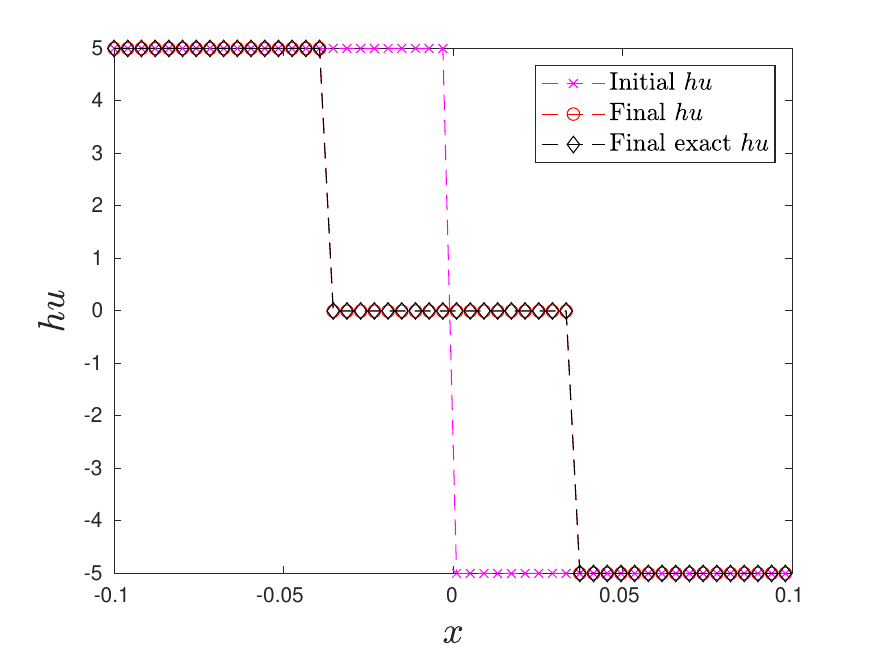}
		\hspace*{1mm}\includegraphics[height=3.8cm, keepaspectratio]{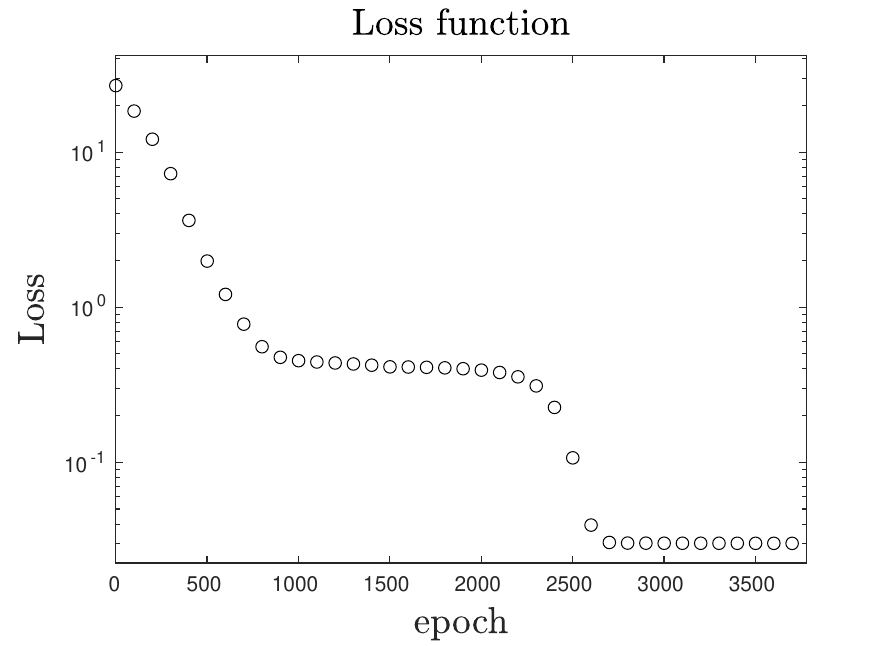}
		\caption{{\bf Experiment 8b.} (Left) Approximate component $x \mapsto  h(x,T)$. (Middle) Approximate component   $x \mapsto hu(x,T)$. (Right) Loss function.}
		\label{hyp_sys_init0b}
	\end{center}
\end{figure}
This experiment shows that the proposed methodology allows for the computation of the solution to (at least simple) Riemann problems.
\paragraph{\bf Experiment 9} In this last experiment, we consider Euler's equations modeling compressible inviscid fluid flows. This is a 3-equation HCL which reads as follows (in conservative form)
\begin{eqnarray*}
	\left.
	\begin{array}{lcl}
		\partial_t \rho + \partial_x(\rho u) & = & 0,\\
		\partial_t (\rho u) + \partial_x(\rho u^2 + P) & = & 0,\\
		\partial_t (\rho E) + \partial_x(\rho Eu  + Pu) & = & 0,
	\end{array}
	\right.
\end{eqnarray*}
with the equation of states $P=\rho(\gamma-1)(E-u^2/2)$ (perfect gas law), $\gamma=1.4$, and where $\rho$ denotes the fluid density, $u$ denotes the fluid velocity and $E$ denotes the total energy. Dirichlet boundary conditions are imposed at the boundary of the spatial domain $(0,2)$.  In this experiment, we consider a stiff benchmark where the solution is constituted by a 1-shock, 3-shock waves and 2-rarefaction wave \cite{MYV}.  Following the strategy proposed in Subsection \ref{subsec:IWD}, we consider the initial data given by
\begin{eqnarray*}
	(\rho_0,\rho_0u_0,\rho_0E_0)=
	\left\{ 
	\begin{array}{ll}
		(0.01, 21, 25632), & x<x_0,\\
		(0.06, 15,119000), & x_0<x<x_1, \\
		(0.235, 60, 125000), & x_1<x<x_2, \\
		(0.14, 0,55750), & x_2<x<x_3, \\
	\end{array}
	\right.
\end{eqnarray*}
with $x_0=0.95$, $x_1=1.0$, $x_2=1.05$ and $x_3=2$ and $T=0.001$ . The initial density, velocity and pressure, are represented in Fig. \ref{hyp_sys_init}.
\begin{figure}[!ht]
	\begin{center}
		\hspace*{1mm}\includegraphics[height=3.8cm, keepaspectratio]{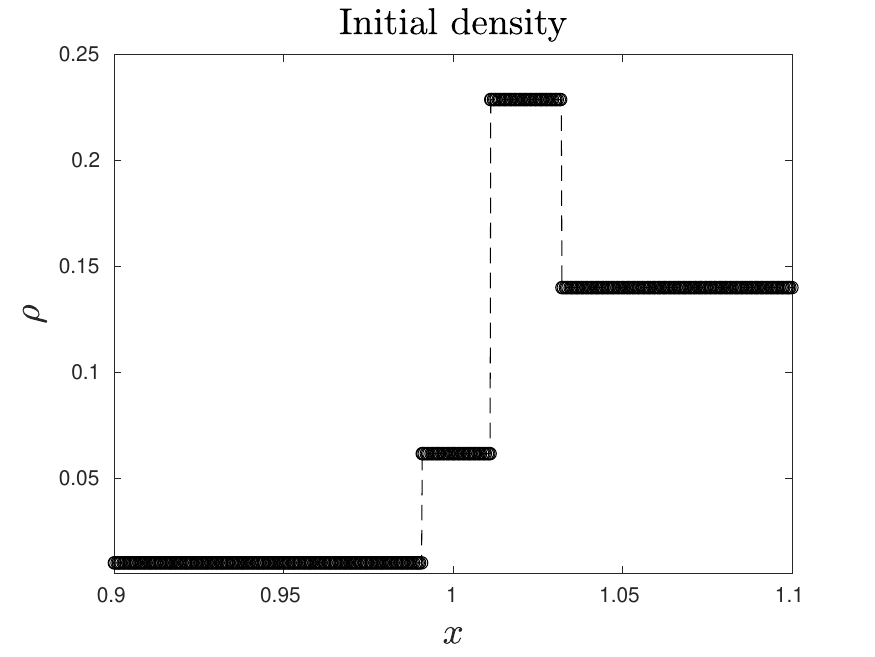}
		\hspace*{1mm}\includegraphics[height=3.8cm, keepaspectratio]{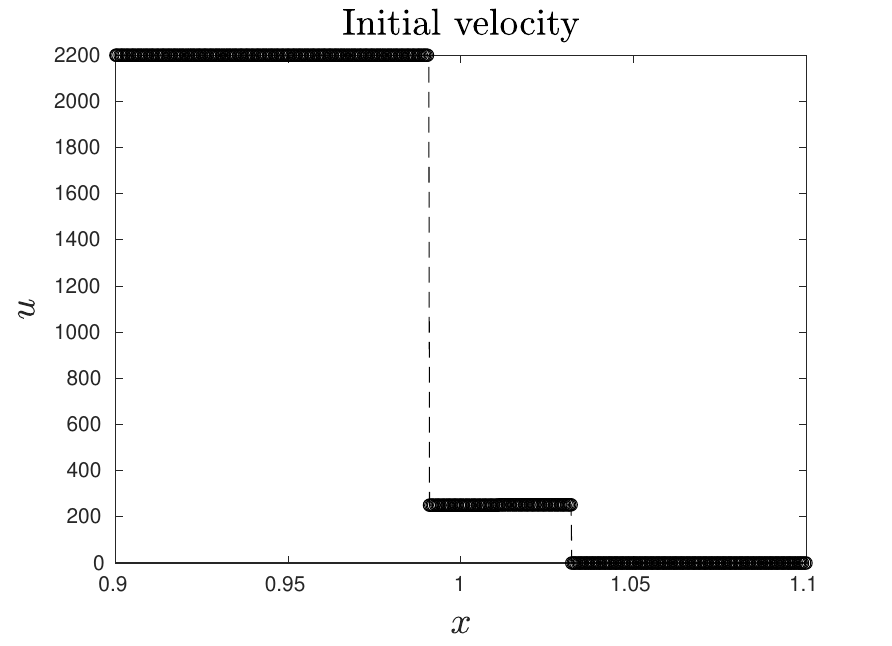}
		\hspace*{1mm}\includegraphics[height=3.8cm, keepaspectratio]{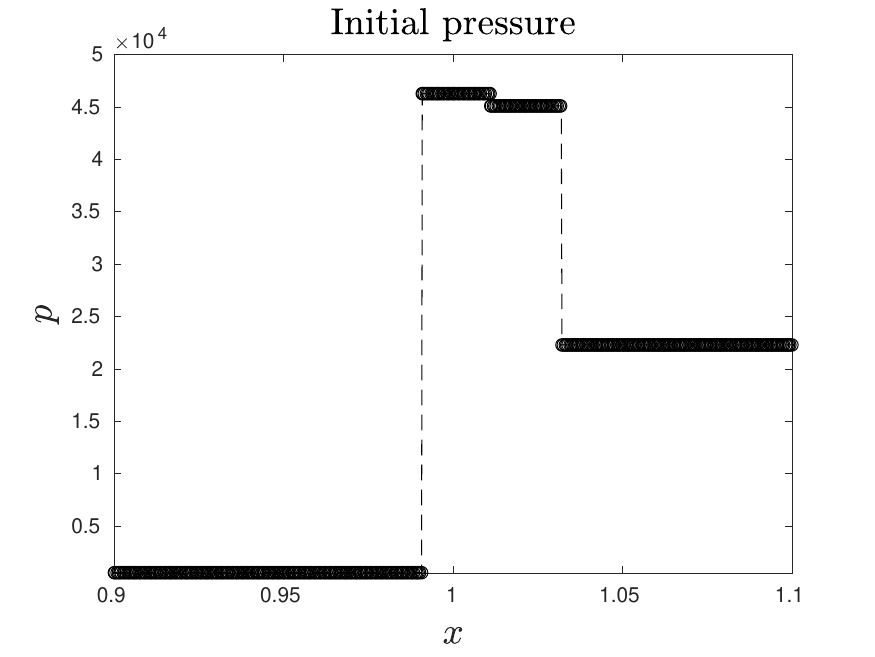}
		\caption{{\bf Experiment 9.} Initial data (Left)  Density. (Middle) Velocity (Right) Pressure.}
		\label{hyp_sys_init}
	\end{center}
\end{figure}

We implement the method developed in Subsection \ref{sec:nd-solver} with $m=3$, $1$ hidden layer and $30$ neurons for each conservative component $\rho$, $\rho u$, $\rho E$ and for the $3$ lines of discontinuity. In Fig. \ref{hyp_sys} we report  the density, velocity and pressure at initial and final times $T$. This test illustrates the precision of the proposed approach, with in particular an accurate approximation of the 2-contact discontinuity which is often hard to obtain with standard solvers.
\begin{figure}[!ht]
	\begin{center}
		\hspace*{1mm}\includegraphics[height=3.8cm, keepaspectratio]{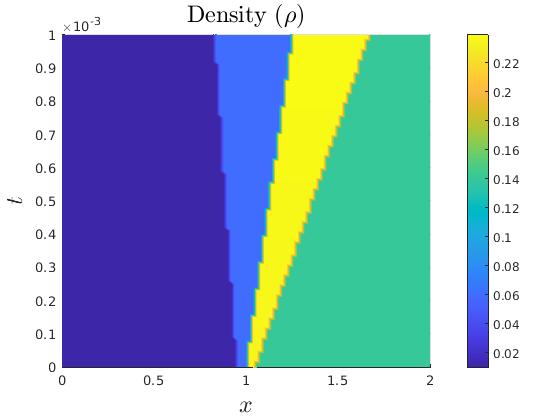}
		\hspace*{1mm}\includegraphics[height=3.8cm, keepaspectratio]{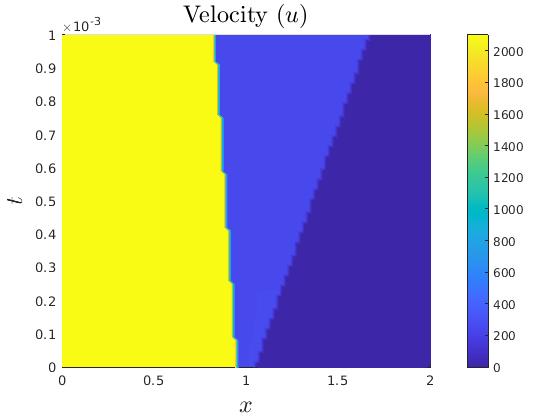} 
		\hspace*{1mm}\includegraphics[height=3.8cm, keepaspectratio]{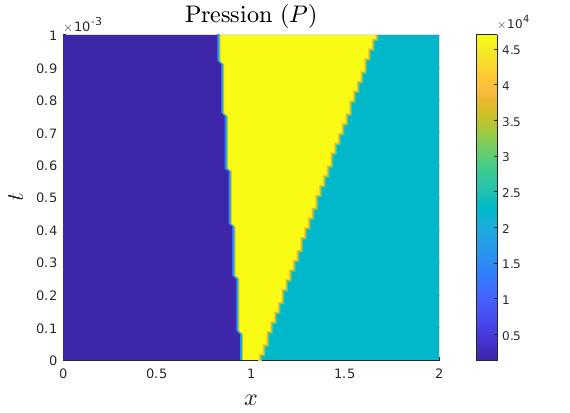}
		\caption{{\bf Experiment 9.} Solution in $(0,2)\times (0,0.001)$. (Left)  Density. (Middle) Velocity (Right) Pressure.}
		\label{hyp_sys}
	\end{center}
\end{figure}
% end

\section{Conclusion}\label{sec:conclusion}
In this paper, we have proposed an original method for solving hyperbolic conservation laws using a non-diffusive neural network method.  The principle of the method is to track the DLs and solve the conservation laws in the subdomains the DLs define, and where the solution is smooth. We have used neural networks to approximate the DLs and the solution of conservation laws in each of the subdomains. The networks are trained by minimizing a loss functional that measures the (norm of the) residues of conservation laws, boundary and 
initial conditions and Rankine-Hugoniot conditions. This approach allows for a computation of shock waves without using a weak formulation of HCL. Other functions could be used to approximate the solution and the DLs, but neural networks provide 
interesting features. Indeed, they allow for automatic differentiation, which avoids the approximation errors for the derivatives, facilitates the implementation of the algorithm,
and allows for an accurate$/$diffusion-free computation of shock waves, solutions to nonlinear hyperbolic conservation laws. 

When the global loss functional approach is slow to converge, a rapidly convergent and embarrassingly parallel 
(Schwarz) domain decomposition method can be used. The latter allows for a decoupling of the optimization procedure thanks to the optimization of local neural networks approximating from which the global approximate solution is reconstructed, as shown in \cite{HypSWR2}.

In a future work, we plan to apply to this methodology to higher dimensional problems.

\section*{CRediT authorship contribution statement}
The authors have contributed equally.

\section*{Declaration of competing interest}
The authors declare that they have no known competing financial interests or personal relationships that could have
appeared to influence the work reported in this paper.

\section*{Data availability}
No data was used for the research described in this article.

%\section*{References}
%\begin{thebibliography}
\bibliographystyle{unsrt}
\nocite{*}

\section*{\refname}

%\bibliography{refs_revision}
\bibliography{jcp-2}
%\end{thebibliography}
 	
\end{document}